\documentclass[preprint,3p,times,numbers]{elsarticle}

\usepackage{amssymb}
\usepackage{amsmath}
\usepackage{enumitem}
\setlist[itemize]{noitemsep,nolistsep}
\setlist[enumerate]{noitemsep,nolistsep}
\usepackage{amsthm}
\usepackage[normalem]{ulem}
\usepackage[utf8]{inputenc}
\usepackage[T1]{fontenc}
\usepackage[nottoc]{tocbibind}

\usepackage{mathtools}
\usepackage{dsfont}
\usepackage{mathrsfs}
\usepackage{csquotes}

\usepackage{hyperref} 
\hypersetup{
  breaklinks=true,
  colorlinks=true,
  linkcolor=blue,
  urlcolor=orange!75!black,
  citecolor=green,
  hypertexnames=false,
  }
\usepackage[noabbrev,capitalize,nameinlink]{cleveref}
\usepackage{thmtools}

\newtheorem{theorem}{Theorem}[section]
\newtheorem{proposition}[theorem]{Proposition}
\newtheorem{corollary}[theorem]{Corollary}
\newtheorem{lemma}[theorem]{Lemma}
\theoremstyle{definition}
\newtheorem{definition}[theorem]{Definition}
\newtheorem{remark}[theorem]{Remark}
\newtheorem{assumption}[theorem]{Assumption}
\newtheorem{example}[theorem]{Example}
\newtheorem{cexample}[theorem]{Counterexample}

\newcommand*\dif{\mathop{}\!\mathrm{d}}
\newcommand{\cL}{\mathcal{L}^0}
\newcommand{\cLs}{\mathcal{L}^0_{\mathrm{s}}}
\newcommand{\cLbs}{\mathcal{L}^0_{\mathrm{bs}}}
\newcommand{\overlinecL}{\overline{\mathcal{L}}{}^0}
\newcommand{\overlinecLs}{\overline{\mathcal{L}}{}^0_{\mathrm{s}}}
\newcommand{\overlinecLbs}{\overline{\mathcal{L}}{}^0_{\mathrm{bs}}}
\newcommand{\cLb}{\mathcal{L}^0_{\mathrm{b}}}
\newcommand{\Lb}{L^0_{\mathrm{b}}}
\newcommand{\cLph}{\mathcal{L}^p_h}
\newcommand{\overlinecLph}{\overline{\mathcal{L}}{}^p_h}
\newcommand{\Lph}{L^p_h}
\newcommand{\Ccr}{\mathcal{C}_{\mathrm{cr}}}
\newcommand{\Cc}{\mathcal{C}_{\mathrm{c}}}
\newcommand{\Crc}{\mathcal{C}^r_{\mathrm{c}}}
\newcommand{\Crcr}{\mathcal{C}^r_{\mathrm{cr}}}

\DeclareMathOperator*{\esssup}{ess\,sup}

\makeatletter 
\def\ps@pprintTitle{
 \let\@oddhead\@empty
 \let\@evenhead\@empty
 \def\@oddfoot{}
 \let\@evenfoot\@oddfoot}
\makeatother

\begin{document}

\begin{frontmatter}

\title{Nonlinear Lebesgue spaces: Dense subspaces, completeness and separability}

\corref{cor}
\cortext[cor]{Corresponding author}

\affiliation[1]{organization={Université Paris Cité, CNRS, MAP5},
 postcode={F-75006},
 city={Paris},
 country={France}}
 \affiliation[2]{organization={Centre Giovanni Borelli, \'Ecole Normale Supérieure Paris-Saclay, Université Paris-Saclay},
 postcode={F-91190},
 city={Gif-sur-Yvette},
 country={France}}

\author[1]{Guillaume Sérieys}
\ead{guillaume.serieys@u-pariscite.fr}

\author[2]{Alain Trouvé}
\ead{alain.trouve@ens-paris-saclay.fr}

\begin{abstract}
$L^p$ spaces of mappings taking values in arbitrary metric spaces, which we call \emph{nonlinear Lebesgue spaces}, play an important role in several fields of mathematics. For instance, membership in these spaces is typically required for transport maps in optimal transport theory and for stochastic processes in probability theory. Nonlinear Lebesgue spaces also arise naturally in applications such as medical imaging, where the physical signals at play often exhibit little regularity and take their values in nonlinear spaces. Yet, these spaces remain little studied in the literature, likely due to their lack of differential structure outside the case where mappings are valued in a linear space. This paper is the first in a series devoted to the study of geometric and analytic properties of nonlinear Lebesgue spaces. The present article exposes a systematic treatment of their measure-theoretic properties, unifying and refining scattered results from the literature while also extending classical results from the linear setting to this broader nonlinear framework---including the characterizations of their completeness and their separability as well as the density of some of their subspaces: the spaces of simple, continuous and smooth mappings. 
\end{abstract}


\end{frontmatter}

\tableofcontents

\section{Introduction}
\label{sec:intro}

Physical phenomena often involve signals with values in nonlinear spaces. For instance, in medical imaging, pixel values may belong to more complex spaces than the real line such as the space of symmetric positive definite matrices as is the case in diffusion tensor imaging \cite{basser1994mr} (see also \cite[Chapter~3]{pennec2019riemannian} for a complete exposition on tensor-valued image processing). Although this space remains convex when equipped with the Euclidean metric, it is often more desirable to equip it with more \emph{meaningful} Riemannian metrics, which can provide desirable geometric properties such as geodesic completeness and nonpositive curvature, and avoid some of the limitations of Euclidean interpolation and extrapolation (see \cite{thanwerdas2023n} for a further discussion on this matter). Another good example is when the space of values is the probability simplex such, as in the case of soft segmentation maps, also called probability maps, to label different tissues while accounting for uncertainty \cite{ashburner2005unified} or to construct probabilistic atlases of organs \cite{mazziotta1995probabilistic}. The probability simplex is usually equipped with the Fisher--Rao metric, which makes it positively curved (see \cite[Section 2.1]{aastrom2017image} for a complete exposition on the geometry of the probability simplex). Outside the manifold-valued case, an example is that of the space of probability measures, which does not carry a differential structure, encountered in Q-ball imaging \cite{tuch2004q} or constrained spherical deconvolution \cite{tournier2004direct} (see also \cite{vogt2018measure} for a complete exposition on the interest of considering measure-valued mappings). Despite the lack of differential structure of this space, it can still be equipped with a metric (in the metric space sense) such as the Hellinger--Kakutani \cite{hellinger1909neue,kakutani1948equivalence} or Kantorovich--Wasserstein \cite{kantorovich1942translocation,villani2009optimal} metrics. Furthermore, medical images cannot usually be assumed to be continuous due to the multiphase nature of anatomy or the inherent non-continuous nature of physical signals. In the case of natural images, the total variation is known to blow up to infinity as the resolution of images increases \cite{gousseau2001natural}, which even rules out the space of mappings with bounded variation to model the space of natural images. All examples mentioned previously, where the domain of definition of the mappings is flat, eventually generalize to cases where the domain of definition is curved, for instance in the framework of functional shapes \cite[Definition~2.1.1]{charon2013analysis}, or in surface-based neuroimaging, where cortical data can be represented as signals defined on spherical parameterizations of cortical surfaces \cite{fischl1999cortical}. For all those reasons, studying merely measurable mappings defined on measurable spaces with values in arbitrary metric spaces is of importance from an applicative viewpoint. 

$L^p$ spaces of mappings taking values in arbitrary metric spaces also appear in more theoretical contexts. In \cite[Section~1.1]{korevaar1993sobolev}, N.~J. Korevaar and R.~M. Schoen introduced such spaces in the study of Sobolev spaces of mappings taking values in arbitrary metric spaces, and also mentioned sufficient conditions for their completeness. This work was followed by several contributions including the work of K.~T.~Sturm (see \cite[Section~3]{Sturm2001} and \cite[Section~4]{sturm2002nonlinear}) and J.~Jost (see \cite[p.~188]{jost1994equilibrium}, \cite[Section~2]{jost1997generalized}  and \cite[Section~4.1]{jost1997nonpositive}), which also address the matter of geodesics and curvature when mappings are valued in a space of non-positive curvature in the sense of Alexandrov (see \cite[Chapter 4]{burago2022course} for a definition). These results are gathered by M.~Bačák in \cite[pp.~18-19]{bacak2014convex}. Such $L^p$ spaces are also used in the literature of optimal transport theory to model spaces of transport maps (see the book by L.~Ambrosio, N.~Gigli and G.~Savaré \cite[Section~5.4]{ambrosio2005gradient}), that is, using the probability theory terminology, spaces of nonlinear random variables with finite $p$-th moment. They also appear as metric completions of certain $L^p$-type geometries. A prominent example is provided by the $L^2$ metric on the manifold of smooth Riemannian metrics, also known as the Ebin metric: its metric completion was characterized by B.~Clarke \cite{clarke2009thesis,clarke2013completion}. More recently, N.~Cavallucci gave a shorter proof of this result, exhibiting the completion explicitly as an $L^2$ space of mappings taking values in a fixed $\operatorname{CAT}(0)$ space \cite{cavallucci2023l2}. A similar completion result was obtained by N.~Cavallucci and Z.~Su for the space of full-rank vector-valued one-forms equipped with a generalized Ebin metric \cite{cavallucci2023metric}. Related recent developments include the studies by D.~Lenze of isometric rigidity of $L^2$ spaces of manifold-valued mappings \cite{lenze2025isometric} and of isometries of the Ebin metric \cite{lenze2025isometries}. These $L^p$ spaces of mappings taking values in arbitrary metric spaces are also useful as theoretical
tools  for the study of curves in Wasserstein spaces by S.~Lisini (see \cite[Section~2.3]{lisini2007characterization} or \cite[Section~1.2]{lisini2006absolutely}). More recently, M.~Bauer, F.~Mémoli, T.~Needham and M.~Nishino \cite[Section~2.1.3]{bauer2025z} derived additional properties of these spaces, of which separability is particularly noteworthy, as, to the best of our knowledge, it had not appeared previously in the literature. On a more applied note, such spaces were used by A.~Effland, S.~Neumayer, J.~Persch, M.~Rumpf and G.~Steidl \cite{neumayer2018,effland2020convergence} as part of S.~Neumayer's PhD thesis \cite{neumayer2021deformation} to model manifold-valued images in the case where the space of values is a Cartan--Hadamard manifold (see \cite[Section~12,~p.~352]{lee2018introduction} for a definition).
Let us conclude with a note on the terminology. As highlighted in \cite[TVS V.86]{bourbaki2013topological}, $L^p$ spaces of real-valued mappings were first introduced by F.~Riesz in \cite{riesz1910untersuchungen}. Yet, they are commonly referred to as \enquote{Lebesgue} spaces \cite[Chapter III, Section 3]{dunford1988linear} or \enquote{Bochner} spaces when mappings are valued in a Banach space \cite[Section~1.2.b]{hytonen2016analysis}. For consistency with the literature, we use the terminology of \enquote{nonlinear Lebesgue spaces} to designate $L^p$ spaces of mappings taking values in arbitrary metric spaces, which was first used by M.~Bačák in \cite[p.~18]{bacak2014convex}.

This article is the first in a series devoted to the study of geometric and analytic properties of nonlinear Lebesgue spaces, with the second article \cite{serieys2026nonlinear} focusing on their curves and geometric properties. In the present article, we aim to gather several measure-theoretic properties of these spaces while generalizing properties usually derived in the \emph{linear case}, that is, when mappings are valued in a linear space. References for the linear case are for instance the books by H.~Brézis \cite{brezis2011functional}, D.~L.~Cohn \cite{cohn2013measure} or T.~Hytönen, J.~van Neerven, M.~Veraar and L.~Weis \cite{hytonen2016analysis}. In particular, this article gathers known facts about these spaces that are spread throughout the literature and proposes a unified exposition while trying as much as possible to rely on minimal assumptions. Most results are thus, to the best of the knowledge of the authors, not present in the literature at this level of generality. For the sake of clarity, the novelty of most results, especially regarding the generality of the assumptions, is therefore discussed in light of existing particular cases from the literature. 

The contributions of this work are twofold. (i) First, it provides necessary and sufficient conditions
for nonlinear Lebesgue spaces to be complete or separable. Outside a
degenerate measure-theoretic situation, completeness is shown to be
equivalent to completeness of the \emph{target space}, that is, the space in which mappings take their values (\cref{th:charact_completeness_lebesgue}). Separability is characterized in
terms of the separability of the target space and of measure-theoretic
properties of the \emph{base space}, that is, the space on which mappings are defined (\cref{th:charact_separability_lebesgue_general}). These results both recover and
extend statements previously established under more restrictive
assumptions, mostly in the linear case. (ii) Second, it provides a systematic study of dense subspaces in nonlinear Lebesgue spaces.
The density of simple mappings is first established (\cref{th:density_simple}). The density of continuous
mappings is then proved under suitable topological and measure-theoretic assumptions (\cref{th:density_continuous}). Finally, the density of smooth mappings is established in the setting of Banach manifolds by using smooth partitions of unity as a substitute for the mollification technique used in the linear case (\cref{th:density_smooth}). Counterexamples are also provided to clarify the role of the assumptions appearing in these
statements.

The paper is structured as follows. \cref{sec:nonlinear_lebesgue_spaces} introduces the basic setting of this article regarding minimal assumptions and introduces the core definitions and results needed to construct nonlinear Lebesgue spaces. \cref{sec:fundamental_properties} gathers elementary and more advanced properties, which provide the core arguments for the proofs of the following sections. \cref{sec:characterization} is devoted to the completeness and separability of nonlinear Lebesgue spaces. Finally, \cref{sec:density_simple,sec:density_continuous,sec:density_smooth} study the density of simple and countably valued, continuous, and smooth mappings, respectively. 

\section{Nonlinear Lebesgue spaces}
\label{sec:nonlinear_lebesgue_spaces}
\subsection{Basic assumptions}
\label{sec:basic_setting}

Let us first introduce the basic assumptions of this article.

\begin{assumption}[Basic assumptions]
\label{assum:minimal}
    We consider mappings from a nonempty set $M$, called the \emph{base space}, to a nonempty set $N$, called the \emph{target space}. $M$ and $N$ are assumed to carry the following structures: 
    \begin{enumerate}[label=(\roman*)]
        \item \label{assum:base_space} $(M, \Sigma_M,\mu_M)$ is a measure space, that is, the set $M$ is paired with a $\sigma$-algebra $\Sigma_M$ on $M$ and a measure $\mu_M: \Sigma_M \to [0,\infty]$ on $(M,\Sigma_M)$.
    \item \label{assum:target_space} $(N,d_N)$ is a metric space with finite metric $d_N: N^2 \to [0,\infty)$.
    \end{enumerate}
\end{assumption}

\begin{remark}[On the terminology of \enquote{metric}]
    A metric, in the sense of \ref{assum:target_space} of \cref{assum:minimal}, is sometimes called a \enquote{distance function} \cite[Section~III.1]{dieudonne1960treatise}. Yet, we stick with the word \enquote{metric} from now on. 
\end{remark}

These basic assumptions will not be repeated in the statements of subsequent results; yet, any additional assumption relative to \cref{assum:minimal} will be specified.

Notation and terminology are introduced progressively and collected for reference in \cref{sec:notations}. 

\subsection{Measurable mappings}
\label{sec:measurable_mappings}

Most known measure-theoretic facts in this section (and in this article) are selected from the monographs of D.~L.~Cohn \cite{cohn2013measure}, V.~Bogachev \cite{bogachev2007measure} and D.~H.~Fremlin \cite{fremlin2000measure,fremlin2001measure}. For a similar exposition in the linear case, see \cite[Chapter~24]{fremlin2001measure}.

In this section, we introduce several classes of mappings that will serve as core ingredients in the construction of nonlinear Lebesgue spaces, namely, the set of measurable mappings and its variants. Hereafter, a set $Z\subset M$ is called \emph{$\mu_M$-null} if it is contained in a measurable set of $\mu_M$-measure zero and we denote the collection of all such sets by $\mathcal{Z}_{\mu_M}$.

\begin{definition}[Measurable mappings]
\label{def:measurable_map}
Define:
\begin{enumerate}[label=(\roman*)]
    \item $\cL(M,N)$ the set of \emph{measurable mappings}, that is, all mappings $f: M\to N$ such that for every set $B$ in $\mathcal{B}(N)$, the Borel $\sigma$-algebra of $N$, it satisfies $f^{-1}(B)\in \Sigma_M$.
    \item $\cLs(M,N)$ the set of \emph{separably valued measurable mappings}, that is, all mappings $f\in \cL(M,N)$ such that its range, denoted by $f(M)$, is separable. 
    \item $\cLbs(M,N)$ the set of \emph{separably valued and bounded measurable mappings}, that is, all mappings $f\in \cLs(M,N)$ such that $f(M)$ is bounded. 
    \item $\overlinecL(M,N)$ the set of \emph{$\mu_M$-measurable mappings}, that is, all mappings $f: M\to N$ such that for all $B\in \mathcal{B}(N)$ the preimage $f^{-1}(B)$ belongs to $\overline{\Sigma}_M$, the completion of the $\sigma$-algebra $\Sigma_M$ under $\mu_M$.
    \item $\overlinecLs(M,N)$ the set of \emph{$\mu_M$-essentially separably valued $\mu_M$-measurable mappings}, that is, all mappings $f\in \overlinecL(M,N)$ such that there exists a set $Z$ in $\mathcal{Z}_{\mu_M}$ for which $f(M\setminus Z)$  is separable. 
    \item $\overlinecLbs(M,N)$ the set of \emph{$\mu_M$-essentially separably valued and bounded $\mu_M$-measurable mappings}, that is, all mappings $f\in \overlinecL(M,N)$ such that there exists $Z\in \mathcal{Z}_{\mu_M}$ for which $f(M\setminus Z)$ is separable and bounded.
\end{enumerate}

\end{definition}

\begin{remark}[On the inclusion order of these sets]
\label{rem:inclusion_measurable}
     Note that we have the inclusion $\cL(M,N)\subset \overlinecL(M,N)$, with the reverse inclusion when $\mu_M$ is a complete measure, that is, when every $\mu_M$-null set belongs to $\Sigma_M$, so that $\Sigma_M$ and  $\mu_M$ coincide with their completions $\overline{\Sigma}_M$ and $\bar{\mu}_M$. Similarly, we have $\cLs(M,N)\subset \overlinecLs(M,N)$ and $\cLbs(M,N) \subset \overlinecLbs(M,N)$. In addition, if $\mu_M$ is complete and $N$ is separable, then $\overlinecLs(M,N)\subset \cLs(M,N)$. Also, when $N$ is separable, $\cLs(M,N)= \cL(M,N)$ and $\overlinecLs(M,N)= \overlinecL(M,N)$.
\end{remark}

\begin{remark}[\enquote{Borel measurable} mappings]
When $M$ is a topological space with $\Sigma_M = \mathcal{B}(M)$, the Borel $\sigma$-algebra of $M$, a mapping belonging to $\cL(M,N)$ is usually called \enquote{Borel measurable}.
\end{remark}

A first well-known result on this class of mappings is that the pointwise limit of measurable and separably valued mappings, when it exists, is itself measurable and separably valued. 

\begin{proposition}[$\cLs(M,N)$ is closed under pointwise limit]
\label{prop:measurability_pointwise_limit}
    Let $(f_n)_{n\in\mathbb{N}}$ be a sequence in $\cLs(M,N)$ and assume that $f(x)\coloneqq \lim_{n\to\infty} f_n(x)$ exists in $N$ for all $x\in M$. Then, $f\in \cLs(M,N)$.
\end{proposition}
\begin{proof}
    See \ref{appendix:measurability_pointwise_limit}.
\end{proof}

In particular, \cref{prop:measurability_pointwise_limit} holds for sequences of \emph{simple mappings}, that is, measurable mappings with a finite range.

\begin{definition}[Simple mappings]
    \label{def:simple_map}
    Define:
\begin{enumerate}[label=(\roman*)]
    \item $\mathcal{E}(M,N)$ the set of \emph{simple mappings}, that is, all mappings $f\in \cL(M,N)$ such that $\lvert f(M)\rvert$, the cardinality of the range of $f$, is finite.
    \item $\overline{\mathcal{E}}(M,N)$ the set of \emph{$\mu_M$-essentially simple mappings}, that is, all mappings $f\in \overlinecL(M,N)$ such that there exists $Z\in \mathcal{Z}_{\mu_M}$ for which $\lvert f(M\setminus Z)\rvert < \infty$.
\end{enumerate}

Note that both $\mathcal{E}(M,N)\subset\cLbs(M,N) $ and $\overline{\mathcal{E}}(M,N)\subset\overlinecLbs(M,N) $.
\end{definition}

Simple mappings are sometimes equivalently defined through a partition of the base space. 

\begin{proposition}[Equivalent definition of simple mappings]
\label{prop:equivalent_def_simple}
    A mapping $f: M\to N$ belongs to $\mathcal{E}(M,N)$ if and only if there exists a finite partition $(M_i)_{i\in I}$ of $M$ such that $M_i \in \Sigma_M$ and $f|_{M_i} \equiv y_i \in N$.
\end{proposition}

\begin{proof} 
See \ref{appendix:equivalent_def_simple}.
\end{proof}

Showing that a mapping belongs to $\cLs(M,N)$ thus usually comes down to verifying that it is the pointwise limit of a sequence of simple mappings.

Now, one can define an equivalence relation on the set of $\mu_M$-measurable mappings.
\begin{proposition}[Equivalence relation on $\overlinecL(M,N)$]
\label{prop:equivalence}
    Define the relation $\sim$ on $\overlinecL(M,N)$ by  $$f\sim f' \iff f(x)=f'(x)\quad\text{for $\mu_M$-a.e.~$x\in M$.}$$ Then, $\sim$ is an equivalence relation on $\overlinecL(M,N)$.
\end{proposition}

Using this equivalence relation, one can define the set of equivalence classes of measurable mappings and its variants.

\begin{definition}[Equivalence classes of measurable mappings]
\label{def:equiv_class_measurable}
    Let
    $$q: \overlinecL(M,N)\to \overlinecL(M,N)/\sim$$ 
    be the quotient mapping associated to $\sim$, given by 
    $$q(f)\coloneqq [f], \quad [f]\coloneqq \{ f'\in \overlinecL(M,N) : f\sim f'\}.$$
    Then, define:
    \begin{enumerate}[label=(\roman*)]
        \item $L^0(M,N) \coloneqq q(\cLs(M,N))$ the set of equivalence classes admitting a separably valued measurable representative.
        \item $\Lb(M,N)\coloneqq q(\cLbs(M,N))$ the set of equivalence classes admitting a separably valued and bounded measurable representative.
        \item $E(M,N) \coloneqq q(\mathcal{E}(M,N))$ the set of equivalence classes admitting a simple representative.
    \end{enumerate}
    
\end{definition}

\begin{remark}[On the separably valued assumption]
\label{rem:separably_valued}
    The assumption that the mappings are separably valued ensures that $d_N$ is Borel measurable when restricted to the separable closed (hence Borel measurable) set $\overline{f(M)}\times \overline{f'(M)}$ for all $(f,f')\in \cLs(M,N)^2$. This allows avoiding Nedoma's pathology cases (see \cite[Sections~15.9 and 15.10]{schilling2021counterexamples} for a complete exposition on this matter) and ensuring that $x\mapsto d_N(f(x),f'(x))$ is always measurable, for any such pair of mappings. In addition, the assumption of separable range is coherent with Bochner integration theory (see \cite[Section~1]{hytonen2016analysis} for a complete exposition on the topic) in which mappings are assumed separably valued to ensure that strong measurability (see \cite[Definition~1.1.4]{hytonen2016analysis} for a definition) and measurability in the sense of \cref{def:measurable_map} coincide (see \cite[Corollary~1.1.10]{hytonen2016analysis} and \cite[Theorem~1.1.20]{hytonen2016analysis}). 
\end{remark}

Another standard result is the fact that a $\mu_M$-measurable mapping is $\mu_M$-essentially separably valued if and only if there exists an equivalent measurable and separably valued mapping (see \cite[Exercise~1.3]{ambrosio2000functions} in the case of a separable target space). 

\begin{proposition}[$\cLs(M,N)$ representative of $\overlinecLs(M,N)$ mappings]
\label{prop:measurable_representative}
    A $\mu_M$-measurable mapping $f$ belongs to $\overlinecLs(M,N)$ if and only if there exists $\Tilde{f}$ in $\cLs(M,N)$ such that $ \tilde{f}\sim f$.
\end{proposition}
\begin{proof}
    See \ref{appendix:proof_measurable_representative}.
\end{proof}

\cref{prop:measurable_representative} will thus play an important role in the definition of nonlinear Lebesgue spaces since it essentially states that up to choosing an appropriate representative in the equivalence class, we can always work with measurable and separably valued mappings. A direct consequence of \cref{prop:measurable_representative} is that the sets introduced in \cref{def:equiv_class_measurable} can equivalently be defined using $\mu_M$-measurable mappings.

\begin{corollary}[Equivalent definitions of the spaces in \cref{def:equiv_class_measurable}]
\label{prop:equiv_def_measurable_class}
    We have:
    \begin{enumerate}[label=(\roman*)]
        \item \label{itm:l0_equiv_def} $L^0(M,N) = q(\overlinecLs(M,N))$
        \item \label{itm:lb_equiv_def} $\Lb(M,N) = q(\overlinecLbs(M,N))$
        \item \label{itm:e_equiv_def} $E(M,N) = q(\overline{\mathcal{E}}(M,N))$
    \end{enumerate}
\end{corollary}

From these first building blocks, we move on to the construction of nonlinear Lebesgue spaces.

\subsection{Lebesgue mappings}
\label{sec:lebesgue_mappings}

Before defining nonlinear Lebesgue spaces, let us first define the $\mathcal{L}^p$ semi-metrics (see \cite[Definition~1.1.4]{burago2022course} for a definition of the notion of semi-metric).

\begin{definition}[$\mathcal{L}^p$ semi-metrics] Let $p\in [1,\infty]$. Then, define:
\begin{enumerate}[label=(\roman*)]
    \item the mapping $D_p: \cLs(M,N)^2\to [0,\infty]$ as 
$D_p(f,f') \coloneqq \lVert d_N(f,f')\rVert_{p,\mu_M}$
where $d_N(f,f')(x)\coloneqq d_N(f(x),f'(x))$. Precisely, for $p\in [1,\infty)$,
$$D_p(f, f') = \left(\int_M d_N(f(x), f'(x))^p\dif \mu_M(x)\right)^{1/p}$$
and, for $p=\infty$, 
$$D_\infty(f,f') = \mu_M\text{-}\esssup_{x\in M} d_N(f(x),f'(x))$$
with the $\mu_M$-essential supremum defined as
$$\mu_M\text{-}\esssup_{x\in M} d_N(f(x),f'(x))\coloneqq \inf\left\{C\in\mathbb{R}_+: d_N(f(x),f'(x)) \leq C \text{ for $\mu_M$-a.e. $x\in M$}\right\}.$$

\item the mapping $\overline{D}_p: \overlinecLs(M,N)^2 \to [0,\infty]$ by replacing $\mu_M$ with $\bar{\mu}_M$ in the definition of $D_p$.
\end{enumerate}
\end{definition}

\begin{remark}[\enquote{$\mathcal{L}^p$ distance} between measurable mappings]
\label{rem:eval_equiv_class}
 For any pair of mappings $(f, f') \in \cLs(M,N)^2$, we call \enquote{$\mathcal{L}^p$ distance} between $f$ and $f'$ the value $D_p(f, f')$. Also, note that, in that case, $D_p(f,f') = \overline{D}_p(f,f')$. Similarly, the \enquote{$\mathcal{L}^p$ distance} between any two equivalence classes $([f],[f'])\in L^0(M,N)$ is given by $D_p([f],[f'])\coloneqq D_p(f,f')$, that is, the evaluation of $D_p$ for two representatives $(f,f')\in \cLs(M,N)^2$. The same holds if only one of the two mappings is an equivalence class.
\end{remark}

At this point, $D_p$ satisfies the symmetry and triangle inequality axioms of a metric on $\cLs(M,N)$ thanks to the symmetry of $d_N$ and the fact that both $d_N$ and $\lVert \cdot \rVert_{p,\mu_M}$ satisfy the triangle inequality for $p\in [1,\infty]$. However, $D_p$ is not necessarily finite and cannot distinguish two distinct measurable mappings that agree $\mu_M$-almost everywhere. In fact, $D_p$ separates equivalence classes in $L^0(M,N)$. 

\begin{proposition}[$\mathcal{L}^p$ semi-metrics separate equivalence classes]
\label{prop:lp_separates_equiv_class}
    Let $p\in [1,\infty]$ and $(f,f')\in \cLs(M,N)^2$. Then, $D_p(f,f')= 0$ if and only if $f\sim f'$.
\end{proposition}

\begin{remark}[$\mathcal{L}^p$ semi-metrics metrize $L^0(M,N)$]
    Note that $L^0(M,N)$ equipped with $D_p$, for any choice of $p\in [1,\infty]$, satisfies all the axioms of a metric space, but with a metric $D_p$ that might be infinite.
\end{remark}

First, to make $D_p$ finite, consider its restriction to \emph{$\mathcal{L}^p$ mappings}, which are defined as mappings at a finite $\mathcal{L}^p$ distance from a \emph{base mapping} $h$ in $L^0(M,N)$.

\begin{definition}[$\mathcal{L}^p$ spaces] Let $h\in L^0(M,N)$ (recall that $h$ has a representative with separable range by \cref{def:equiv_class_measurable}) and $p\in [1,\infty]$. Then, define:
\begin{enumerate}[label=(\roman*)]
    \item $\cLph(M,N)\coloneqq\{ f\in \cLs(M,N): D_p(f, h) < \infty\}$ the set of measurable mappings at a finite $\mathcal{L}^p$ distance from $h$ (\cref{rem:eval_equiv_class}). 
\item $\overlinecLph(M,N)$ using $\overline{D}_p$ instead of $D_p$ in the definition of $\cLph(M,N)$.
\end{enumerate}

\end{definition}

\begin{remark}[Choice of the base mapping in the linear case]
    When $N$ is a normed vector space, the base mapping $h$ is usually taken as $h\equiv 0_N$, with $0_N$ the identity element of $N$. 
\end{remark}

\begin{remark}[On the inclusion order of $\mathcal{L}^p$ spaces]
\label{rem:bounded_measurabe_maps}
    If $h\in \Lb(M,N)$, we have the inclusion $\cLb(M,N) \subset \mathcal{L}^\infty_h(M,N)$ for any choice of measure $\mu_M$. However, $\cLb(M,N) \subset \cLph(M,N)$, $p\in [1,\infty)$, only holds when $\mu_M$ is finite. In addition, $\mathcal{L}^{p'}_h(M,N) \subset \cLph(M,N)$, for $1 \leq p \leq p' \leq \infty$, only holds when $\mu_M$ is finite, by Hölder's inequality \cite[Proposition~3.3.2]{cohn2013measure}.
\end{remark}

Then, $(\cLph(M,N),D_p)$ satisfies all axioms of a metric space, except separation: it is a \emph{semi-metric} space. We will often write $\cLph(M,N)$ instead of $(\cLph(M,N),D_p)$. Now, to make it a metric space, $\mu_M$-a.e.~identical mappings should be identified through the previously defined equivalence relation (\cref{prop:equivalence}).

\begin{definition}[Nonlinear Lebesgue spaces]
\label{def:lebesgue_spaces}
Let $h\in L^0(M,N)$ and $p\in [1,\infty]$. Then, the nonlinear $p$-Lebesgue space $L_{h}^p(M,N)$ is defined as the quotient space 
$$L_{h}^p(M,N) \coloneqq q(\cLph(M,N)),$$
where we recall that $f\sim f'$ if and only if $D_p(f,f')=0$ (\cref{prop:lp_separates_equiv_class}).
Also, we will usually omit $p$ when referring to nonlinear $p$-Lebesgue spaces and call \enquote{Lebesgue mappings} the elements of nonlinear Lebesgue spaces. 
\end{definition}

An immediate consequence of \cref{prop:equiv_def_measurable_class} is that nonlinear Lebesgue spaces can equivalently be defined using $\overline{\mathcal{L}}{}^p$ mappings.

\begin{proposition}[Equivalent definition using $\overline{\mathcal{L}}{}^p$ mappings]
\label{prop:equiv_def_lebesgue}
    Let $h\in L^0(M,N)$ and $p\in [1,\infty]$. Then, we have the equivalent definition of nonlinear Lebesgue spaces as 
    $$L_{h}^p(M,N) = q(\overlinecLph(M,N)).$$
\end{proposition}

\begin{remark}[Equivalent definition using $L^0(M,N)$]
    \cref{prop:equiv_def_lebesgue} is obvious from the strictly equivalent definition: 
    $$\Lph(M,N) = \big\{ f\in L^0(M,N): D_p(f,h) < \infty \big\}.$$
\end{remark}

\begin{remark}[Related results in the literature]
    \cref{prop:equiv_def_lebesgue} is also known for the linear case \cite[Exercise 244X (a)]{fremlin2001measure}.
\end{remark}

Thanks to \cref{prop:equiv_def_lebesgue}, we will exclusively work with measurable and separably valued mappings in the rest of the article. 

The construction of nonlinear Lebesgue spaces developed above eventually yields that nonlinear Lebesgue spaces are metric spaces.

\begin{proposition}[Nonlinear Lebesgue spaces are metric spaces]
\label{prop:lebesgue_metric_space}
Let $p\in [1,\infty]$. When equipped with $D_p$, $\Lph(M,N)$ becomes a metric space with finite metric. 
\end{proposition}

In what follows, we usually write $\Lph(M,N)$ when referring to the metric space $(\Lph(M,N),D_p)$.

\begin{proof}[Proof of \cref{prop:lebesgue_metric_space}]
    Let $p\in [1,\infty]$. See \cref{rem:eval_equiv_class} for the evaluation of $D_p$ on equivalence classes. Then, the symmetry and triangle inequality axioms are respectively inherited from symmetry of $d_N$ and the fact that both $d_N$ and $\lVert \cdot \rVert_{p,\mu_M}$ satisfy the triangle inequality. \cref{prop:lp_separates_equiv_class} then yields the separation axiom. Finiteness follows from the construction of nonlinear Lebesgue spaces and the triangle inequality.
\end{proof}

Before going to the next section, we highlight the fact that we will frequently not distinguish representatives from their associated equivalence class, so that we can write $f\in \Lph(M,N)$ while treating $f$ as a $\mathcal{L}^p$ mapping when it is clear from the context. Also, to ensure that our notation remains consistent with the literature on linear Lebesgue spaces, we denote $\mathcal{L}^p(M,N)\coloneqq \cLph(M,N)$ and $L^p(M,N)\coloneqq \Lph(M,N)$ when $N$ is a normed vector space and the base mapping is set to $h\equiv 0_N$ with $0_N$ the identity element of $N$.

\section{Preliminary results}

\label{sec:fundamental_properties}

In this section, we collect results that require no further assumptions beyond the basic setting or occur multiple times as key arguments in the proofs of the main results. We provide proofs of the results that serve as key arguments in the following sections, and postpone the proofs of auxiliary, nonetheless useful, results to the appendix.

\subsection{First elementary results}

A first useful result on measurable mappings is the fact that any separably valued measurable mapping can be $\mu_M$-essentially approximated by countably valued measurable mappings.

\begin{proposition}[Density of countably valued mappings in $(L^0(M,N),D_\infty)$]
\label{prop:countably_finite_radius}
    The set of countably valued measurable mappings determines a dense subspace of $L^0(M,N)$ for the $\mathcal{L}^\infty$ topology.
\end{proposition}

\begin{proof}
    See \ref{appendix:proof_countably_finite_radius}.
\end{proof}

Back to nonlinear Lebesgue spaces, the definition of $L^p_{h}(M,N)$ can be made independent of $h$ by assuming that $\mu_M$ is a finite measure, and restricting to bounded base mappings. 

\begin{proposition}[Conditions of invariance to base mapping]
\label{prop:indep_def}
    Let $p\in [1,\infty]$. Suppose that $\mu_M$ is finite. Then,  we have for all $(h,h')\in \Lb(M,N)^2$ that $\Lph(M,N) = L^p_{h'}(M,N)$.
\end{proposition}

\begin{proof}
\cref{prop:indep_def} appears without proof in the literature (for example, in \cite[Section~3, p.~326]{sturm2002nonlinear}), so we provide one in \ref{appendix:indep_def}.
\end{proof}

An important matter to know about nonlinear Lebesgue spaces is whether they are reduced to the equivalence class of their base mapping, that is, $\Lph(M,N) =\left\{h\right\}$ for $h\in L^0(M,N)$. In such cases, nonlinear Lebesgue spaces are called \emph{trivial}. Nonlinear Lebesgue spaces are then called \emph{nontrivial} when they are strictly bigger, in the sense of inclusion, than the equivalence class of the base mapping. For $p\in [1,\infty)$, this is the case if and only if $\mu_M$ is not purely infinite, that is, the range of $\mu_M$ is not reduced to $\{0,\infty\}$, and $\lvert N\rvert > 1$. For $p=\infty$, this is the case if and only if $\mu_M$ is not trivial, that is, the range of $\mu_M$ is not reduced to $\{0\}$, and $\lvert N\rvert > 1$. 

\begin{proposition}[Nontrivial nonlinear Lebesgue spaces]
\label{prop:trivial_lebesgue}
    Let $h\in L^0(M,N)$. Then, the following assertions hold:
    \begin{enumerate}[label=(\roman*)]
        \item for every $p\in [1,\infty)$, the space $\Lph(M,N)$ is nontrivial if and only if $\mu_M$ is not purely infinite and $\lvert N\rvert > 1$.
        \item the space $L^\infty_h(M,N)$ is nontrivial if and only if $\mu_M$ is nontrivial and $\lvert N\rvert > 1$.
    \end{enumerate}
\end{proposition}
\begin{proof}
    See \ref{appendix:proof_trivial_lebesgue}.
\end{proof}

Also, when dealing with countably many Lebesgue mappings (e.g.~when dealing with sequences), one can avoid assuming that $\mu_M$ is $\sigma$-finite and/or that $h$ is bounded by observing that Lebesgue mappings differ from the base mapping on a sequence $(B_n)_{n\in\mathbb{N}}$ of sets in $\mathcal{F}_{\mu_M}$, the set of measurable sets with finite $\mu_M$-measure, such that $h$ is bounded on each set of this sequence. 

\begin{proposition}[Properties of sets where Lebesgue mappings differ from the base mapping]
\label{prop:lebesgue_differ_from_base_bounded}
Let $p\in [1,\infty)$ and $h\in L^0(M,N)$. Then, any mapping in $\Lph(M,N)$ differs from $h$ on a $\sigma$-finite measurable set composed of sets on which $h$ is bounded, that is, for all $f\in \Lph(M,N)$ there exists a sequence $(B_n)_{n\in\mathbb{N}}$ of sets in $\mathcal{F}_{\mu_M}$ such that, up to a $\mu_M$-null set, $\left\{x\in M: f(x)\neq h(x)\right\}= \cup_{n\in \mathbb{N}} B_n$ and $h|_{B_n}\in \Lb(B_n,N)$ for all $n\in \mathbb{N}$.
\end{proposition}

\begin{proof}
    See \ref{appendix:proof_lebesgue_differ_from_base_bounded}.
\end{proof}

Finally, we collect a few useful regularity results about some mappings of nonlinear Lebesgue spaces in simpler spaces, which will come in handy in a few proofs from the following sections. 

\begin{proposition}[Lipschitz continuity of restriction and pointwise distance to base mapping]
\label{prop:continuity_restriction}
Let $h\in L^0(M,N)$ and $p\in [1,\infty]$. Then, the mappings:
\begin{enumerate}[label=(\roman*)]
    \item $r_B:f\in \cLph(M,N)\mapsto f|_B\in \cLph(B,N)$ for any choice of $B\in \Sigma_M$ and
    \item $\varphi_h: f\in \cLph(M,N)\mapsto d_N(f,h)\in \mathcal{L}^p(M,\mathbb{R}_+)$
\end{enumerate}
are both $1$-Lipschitz continuous. Also, note that, when restricting the base space to a measurable subset, it is always assumed that we take the restriction of the base mapping $h$ to this measurable subset.
\end{proposition}
\begin{proof}
    See \ref{appendix:continuity_restriction}.
\end{proof}

We now move on to more advanced properties used as key arguments in the following sections. 

\subsection{Key results for the main proofs}

A class of mappings that will play an important role in transporting properties from nonlinear Lebesgue spaces to their target space is the space of equivalence classes of constant mappings between $M$ and $N$.
\begin{definition}[Equivalence classes of constant mappings]
\label{def:constant_mappings}
    Define $\Delta(M,N) \coloneqq q(\left\{f\equiv y: y\in N\right\})$ the set of equivalence classes admitting a constant representative, which is a subset of $L^0(M,N)$.
\end{definition}

In particular, when $\mu_M$ is finite and nontrivial, that is, the range of $\mu_M$ is not reduced to $\{0\}$, the target space is isometric to the set $\Delta(M,N)$, which is a closed subset of nonlinear Lebesgue spaces.

\begin{proposition}[Properties of the subspace of constant mappings]
\label{prop:constant_mappings}
Let $h\in \Lb(M,N)$ and $p\in [1,\infty]$. Suppose that $\mu_M$ is finite and nontrivial. Then, the following assertions hold:
\begin{enumerate}[label=(\roman*)]
    \item \label{cond:constant_closed} $\Delta(M,N)$ is a closed subset of $\Lph(M,N)$.
    \item there is a (scaled) isometry between $N$ and $\Delta(M,N)$ with $\mu_M(M)^{1/p}$ as scaling factor for $p\in [1,\infty)$, and no scaling factor for $p=\infty$.
\end{enumerate}

When $p=\infty$, both assertions hold even when $\mu_M$ is not finite.
\end{proposition}
\begin{proof}

    Let $h\in \cLbs(M,N)$ and $p\in [1,\infty]$.

    \begin{enumerate}[label=(\roman*),wide]
        \item Let $f\in \Delta(M,N)$, so that we can assume that $f\equiv z_0\in N$. Then 
    $$D_p(f,h) \leq \begin{cases}
        \mu_M(M)^{1/p}\,\sup_{x\in M} d_N(z_0,h(x)),&\text{if $p\in [1,\infty)$}\\
        \sup_{x\in M} d_N(z_0, h(x)),&\text{if $p=\infty$}
    \end{cases}.$$
    Thus, since $h$ is bounded, $D_p(f,h)<\infty$ for $p\in [1,\infty)$ because $\mu_M$ is finite, whereas $D_\infty(f,h) < \infty$ without any finiteness assumption on $\mu_M$.
    Now, let $(f_n)_{n\in\mathbb{N}}$ be a sequence in $\Delta(M,N)$ that converges to $f\in \Lph(M,N)$ in $\Lph(M,N)$. By \cite[Proposition~3.1.5]{cohn2013measure}, we can assume, up to the extraction of a subsequence, that there exists $Z\in \mathcal{Z}_{\mu_M}$ such that for all $x\in M\setminus Z$ we have $f_n(x)\to f(x)$ as $n\to \infty$ and $f_n$ is constant on $M\setminus Z$. In particular, pick $x_0\in M\setminus Z$, so that, using the fact that $f_n$ is constant on $M\setminus Z$ for all $n\in \mathbb{N}$, we have for all $x\in M\setminus Z$ that $f_n(x_0)\to f(x)$ as $n\to \infty$ and, by the uniqueness of the limit, $f(x)=f(x_0)$ for all $x\in M\setminus Z$, that is, $f\in \Delta(M,N)$.
    \item Let $(f,f')\in \Delta(M,N)^2$. Then, by definition of $\Delta(M,N)$, there exists a unique pair $(y,y')\in N^2$ such that $f(x) =  y$ and $f'(x)= y'$ for $\mu_M$-a.e. $x\in M$. In addition, 
    \begin{equation*}D_p(f,f') = \begin{cases}
        d_N(y,y')\,\mu_M(M)^{1/p},&\text{if $p\in [1,\infty)$}\\
        d_N(y,y'),&\text{if $p=\infty$}
    \end{cases}.\qedhere\end{equation*}
    \end{enumerate}
\end{proof}

When $\mu_M$ is not purely infinite, a notable consequence of \cref{prop:constant_mappings} is that the target space is isometric to a closed subset of nonlinear Lebesgue spaces.

\begin{corollary}[Closed isometric image of the target space into nonlinear Lebesgue spaces]
\label{cor:closed_embedding_target}
    Let $h\in L^0(M,N)$. Then, the following assertions hold:
    \begin{enumerate}[label=(\roman*)]
        \item for every $p\in [1,\infty)$, if $\mu_M$ is not purely infinite, then $N$ is isometric (up to a scaling factor, as mentioned in \cref{prop:constant_mappings}) to a closed subset of $\Lph(M,N)$.
        \item if $\mu_M$ is nontrivial, then $N$ is isometric to a closed subset of $L^\infty_h(M,N)$.
    \end{enumerate}
    
\end{corollary}
\begin{proof}
    Let $h\in\cLs(M,N)$. Fix $p\in [1,\infty)$ in case (i), and set $p=\infty$ in case (ii). In case (i), since $\mu_M$ is not purely infinite, there exists $B\in\Sigma_M$ such that $\mu_M(B)\in (0,\infty)$. In case (ii), since $\mu_M$ is nontrivial there exists $B\in\Sigma_M$ such that $\mu_M(B) > 0$. Then, let $z_0\in N$ and  $B_n \coloneqq B \cap \left\{x\in M: d_N(z_0, h(x))\leq n\right\}$ for every $n\in\mathbb{N}$, so that $(B_n)_{n\in\mathbb{N}}$ is an increasing sequence of measurable sets such that $B=\cup_{n\in\mathbb{N}} B_n$. Hence, by the countable subadditivity of $\mu_M$, we have that $0 < \mu_M(B) \leq \sum_{n\in\mathbb{N}} \mu_M(B_n)$, so that there exists $n_0\in \mathbb{N}$ such that $\mu_M(B_{n_0}) > 0$ and let $A\coloneqq B_{n_0}$ in the rest of the proof. By construction, $h$ is bounded on the measurable set $A$. In addition, in case (i), $A$ has finite $\mu_M$-measure, since $A\subset B$, whereas in case (ii) only the condition $\mu_M(A) > 0$ is needed. Now, define for all $y\in N$ the mapping $f_y: M\to N$ as $f_y|_A\equiv y$ and $f_y|_{M\setminus A} = h|_{M\setminus A}$. Then, recalling that $h$ is bounded on $A$, that $A$ has finite $\mu_M$-measure in case (i), and that this finiteness assumption is not required by \cref{prop:constant_mappings} when $p=\infty$, \cref{prop:constant_mappings} yields that the mapping
    \begin{align*}
        i : \begin{cases}
            N &\to \Lph(M,N)\\
            y &\mapsto [f_y]
        \end{cases}
    \end{align*}
    is isometric up to the scaling factor $\mu_M(A)^{1/p}$ when $p\in [1,\infty)$, and is isometric without a scaling factor when $p=\infty$. In addition, its range $i(N)$ is closed in $\Lph(M,N)$. Indeed, if $(f_{y_n})_{n\in\mathbb{N}}$ converges to $f$ in $\Lph(M,N)$, then $f|_{M\setminus A}\sim h|_{M\setminus A}$, while \cref{prop:constant_mappings} shows that $f|_A$ is $\mu_M$-almost everywhere constant. Hence, there exists $y\in N$ such that $f(x)=f_y(x)$ for $\mu_M$-a.e.~$x\in M$, so that $f\in i(N)$.
    \end{proof}

Another class of mappings that will play an important role in identifying dense subspaces is the set of \emph{almost simple} mappings, that is, measurable mappings that are simple on a measurable set of finite $\mu_M$-measure and ($\mu_M$-a.e.) equal to the base mapping otherwise. 
\begin{definition}[Almost simple mappings]
\label{def:almost_simple}
Let $h\in L^0(M,N)$. Then, define:
    \begin{enumerate}[label=(\roman*)] 
        \item $\mathcal{E}_h(M,N)$ the set of \emph{almost simple mappings}, that is, all mappings $f\in \cLs(M,N)$ such that there exists $B\in \mathcal{F}_{\mu_M}$, that is, a measurable set of finite $\mu_M$-measure, for which $f|_{M\setminus B}\sim h|_{M\setminus B}$ and $f|_B$ is a simple mapping (\cref{def:simple_map}) from $B$, equipped with the trace $\sigma$-algebra $\Sigma_M|_{B}\coloneqq \left\{B\cap A: A\in \Sigma_M\right\}$, to $N$.
        \item $E_h(M,N) \coloneqq q(\mathcal{E}_h(M,N))$ the set of equivalence classes admitting an almost simple mapping representative.
    \end{enumerate}
\end{definition}

Then, without further assumptions, one can show the density of this set in nonlinear Lebesgue spaces for $p\geq 1$.

\begin{proposition}[Density of almost simple mappings for $p\geq 1$]
\label{prop:density_almost_simple}
Let $h\in L^0(M,N)$ and $p\in [1,\infty)$. Then, $E_h(M,N)\cap \Lph(M,N)$ is a dense subspace of $\Lph(M,N)$.
    
\end{proposition}
\begin{proof}
    Let $h\in \cLs(M,N)$, $p\in [1,\infty)$, $f\in \cLph(M,N)$ and $\varepsilon > 0$. $f$ being separably valued, there exists a countable dense set $Q_f \subset f(M)$. Hence, we can assume that $Q_f = \left\{y_n\in f(M):n\in\mathbb{N}\right\}$.
    
    \par\medskip \noindent \emph{Step 1 (Approximation by a mapping that differs from $h$ on a measurable set of finite $\mu_M$-measure):} Define for all $n\in \mathbb{N}^*$ the mapping
$$f_n(x) \coloneqq \begin{cases}
     f(x),&\text{if $d_N(f(x),h(x))\geq n^{-1}$}\\
     h(x),&\text{otherwise}
\end{cases},$$
so that for all $x\in M$ we have $f_n(x)\to f(x)$ as $n\to \infty$. In addition, $d_N(f_n(x),f(x))^p \leq d_N(h(x),f(x))^p$. Hence, since $f\in \cLph(M,N)$, we have, by Lebesgue's dominated convergence theorem \cite[Theorem~2.4.5]{cohn2013measure}, that $D_p(f_n,f)^p\to 0$ as $n\to \infty$. Therefore, pick $n_0\in \mathbb{N}^*$, such that $D_p(f_{n_0},f)< \varepsilon/3$. Thus, $f_{n_0} \in \cLph(M,N)$ and defining the measurable set $A_{n_0} \coloneqq \left\{x\in M: d_N(f(x),h(x))\geq {n_0}^{-1}\right\}$, we get $\mu_M(A_{n_0})<\infty$.

\par\medskip \noindent \emph{Step 2 (Approximation by a bounded mapping on a measurable set of finite $\mu_M$-measure):} Let $R \coloneqq \varepsilon/(3\mu_M(A_{n_0})^{1/p})$ and define
$$g_n(x) \coloneqq \begin{cases}
    f_{n_0}(x),&\text{if $x\in A_{n_0}$ and $f_{n_0}(x)\in \cup_{k=0}^n B(y_k, R)$}\\
    h(x),&\text{otherwise}
\end{cases}.$$
Thus, we have for all $x\in M$ that, by density of $Q_f$, $g_n(x)\to f_{n_0}(x)$ as $n\to \infty$ and $d_N(g_n(x),f_{n_0}(x))^p \leq d_N(h(x),f_{n_0}(x))^p$. Hence, since $f_{n_0}\in \cLph(M,N)$, we have, by Lebesgue's dominated convergence theorem \cite[Theorem~2.4.5]{cohn2013measure}, that $D_p(g_n,f_{n_0})^p\to 0$ as $n\to \infty$. Therefore, pick $n_1\in \mathbb{N}$, such that $D_p(g_{n_1},f_{n_0}) < \varepsilon/3$ and define the set $B_{n_1} \coloneqq \cup_{k=0}^{n_1} B(y_k, R)$ and the measurable set $A_{n_1} \coloneqq A_{n_0}\cap f_{n_0}^{-1}(B_{n_1})$, which has finite $\mu_M$-measure.

\par\medskip \noindent \emph{Step 3 (Approximation by an almost simple mapping):} Now, define the set $P_0\coloneqq B(y_0,R)$ and for all $1\leq n\leq n_1$ the sets $P_n\coloneqq B(y_n,R)\setminus \cup_{k=0}^{n-1} B(y_k,R)$, so that $(P_n)_{0\leq n\leq n_1}$ forms a partition of $B_{n_1}$. Also, define $\chi: B_{n_1}\to \left\{0,1,\ldots, n_1\right\}$ such that for all $y\in B_{n_1}$ it satisfies $y \in P_{\chi(y)}$. Then, define 
$$g(x)\coloneqq \begin{cases}
    y_{\chi\circ g_{n_1}(x)},&\text{if $x\in A_{n_1}$}\\
    h(x),&\text{if $x\notin A_{n_1}$}
\end{cases}.$$
By construction, $g\in \mathcal{E}_h(M,N)$ and we have for all $x\in M$ that $d_N(g(x),g_{n_1}(x)) \leq \mathds{1}_{A_{n_1}}(x) R$. Hence, $D_p(g,g_{n_1}) < \varepsilon/3$.

\noindent Finally, we get, by the triangle inequality, that 
\begin{equation*}
D_p(f,g) \leq D_p(f,f_{n_0}) + D_p(f_{n_0},g_{n_1}) + D_p(g_{n_1},g) < \varepsilon.\qedhere
\end{equation*}
\end{proof}

\begin{remark}[Stronger conclusion on $g(A_{n_1})$ when $N$ is separable]
\label{rem:range_set_density_simple}
    When $N$ is assumed separable, that is, when there exists a countable dense set $Q\subset N$, $g$ can be constructed such that $g(A_{n_1})\subset Q$ by replacing $Q_f$ with $Q$ in the proof. 
\end{remark}

In the following sections, we will see that, in this general setting, we can retrieve many standard properties and density results from the linear case (see \cite[Section~4]{brezis2011functional}, \cite[Section~3.4]{cohn2013measure} or \cite[Section~1.2.b]{hytonen2016analysis}). 

\section{Characterization of completeness and separability}
\label{sec:characterization}

\subsection{Characterization of completeness}
\label{sec:completeness}

A first known fact about nonlinear Lebesgue spaces is that they inherit completeness from their target space \cite[Proposition~3.3]{Sturm2001}.

\begin{proposition}[Completeness of nonlinear Lebesgue spaces]
\label{prop:completeness_lebesgue}
Let $h\in L^0(M,N)$ and $p\in [1,\infty]$. Suppose that $\Lph(M,N)$ is trivial (\cref{prop:trivial_lebesgue}) or that $N$ is complete. Then, $\Lph(M,N)$ is complete.
\end{proposition}

\cref{prop:completeness_lebesgue} is pretty straightforward as it follows from a direct generalization of Riesz-Fischer's proof (see \cite[Theorem~4.8]{brezis2011functional} or \cite[Theorem~3.4.1]{cohn2013measure}) and is proved in several references in the literature under slightly more restrictive assumptions on the base mapping  (see \cite[Section~1.1]{korevaar1993sobolev}, \cite[Proposition~3.3]{Sturm2001} or, more recently, \cite[Proposition~1.2.18]{bacak2014convex}). Still, we provide a detailed proof, based on \cite[Section~1.1]{korevaar1993sobolev} and \cite[pp.~86-87]{federer2014geometric}, for the sake of completeness.

\begin{proof}[Proof of \cref{prop:completeness_lebesgue}] 

\noindent Let $h\in \cLs(M,N)$ and $p\in [1,\infty]$. When $\Lph(M,N)$ is trivial, that is, $\Lph(M,N)=\left\{[h]\right\}$, it is also complete as it contains a single element. Therefore, assume that $\Lph(M,N)$ is nontrivial and that $N$ is complete. Now, let $z_0\in N$ and $(f_n)_{n\in\mathbb{N}}$ be a Cauchy sequence in $\cLph(M,N)$. Up to the extraction of a subsequence (see the proof of \cite[Theorem~4.8]{brezis2011functional}), we can assume that for all $n\in \mathbb{N}$ the following holds
    $$D_p(f_{n+1}, f_{n}) \leq 2^{-n}.$$

\noindent
If $p=\infty$, by definition of $D_\infty$ and since countably many mappings are involved, we can assume there exists $Z\in \mathcal{Z}_{\mu_M}$ such that for all $x\in M\setminus Z$ and $n\in \mathbb{N}$ we have 
$$d_N(f_{n+1}(x),f_n(x))\leq D_\infty(f_{n+1},f_n)\leq 2^{-n}.$$
Hence, $(f_n(x))_{n\in\mathbb{N}}$ is a Cauchy sequence in the complete metric space $N$, so that we can introduce the mapping $f:M\to N$ defined such that for all $x\in M\setminus Z$ it satisfies $f(x)\coloneqq \lim_{n\to\infty} f_n(x)$ and $f|_Z\equiv z_0$. Then, $f$ belongs to $\cLs(M,N)$, by \cref{prop:measurability_pointwise_limit}. In addition, we have for all $x\in M\setminus Z$ that 
$$d_N(f(x),h(x)) = \lim_{n\to \infty}  d_N(f_n(x),h(x))\leq \sup_{n\in\mathbb{N}} D_\infty (f_n,h) < \infty.$$ 
and 
$d_N(f(x),f_n(x)) = \lim_{k\to \infty}  d_N(f_k(x),f_n(x))\leq \sup_{k\geq n} D_\infty(f_k,f_n) \to 0$ as $n\to \infty$. Thus, $f$ even belongs to $\mathcal{L}^\infty_h(M,N)$ and $f_n\to f$ in $\mathcal{L}^\infty_h(M,N)$ as $n\to \infty$.

\noindent If $p\in [1,\infty)$, define for all $n\in \mathbb{N}$ the mapping $g_n(x) \coloneqq \sum_{k\leq n} d_N(f_k (x), f_{k+1} (x))$. Then, $(g_n)_{n\in\mathbb{N}}$ is an increasing sequence of $[0,\infty]$-valued measurable functions, so that, denoting $g(x)\coloneqq \lim_{n\to \infty} g_n(x)$, we have, by Beppo Levi's theorem \cite[Corollary~2.4.2]{cohn2013measure} and the triangle inequality, 
    $$\lVert g\rVert_p  = \lim_{n\to \infty} \lVert g_n\rVert_p \leq \lim_{n\to \infty}\sum_{k\leq n} D_p(f_{k+1},f_k)\leq   \sum_{n\in\mathbb{N}} 2^{-n}  <\infty.$$
    Therefore, by \cite[Corollary~2.3.14]{cohn2013measure}, there exists $Z\in \mathcal{Z}_{\mu_M}$ such that for all $x\in M\setminus Z$ 
    $$g(x) = \sum_{n\in\mathbb{N}} d_N(f_n(x),f_{n+1}(x)) < \infty.$$
    Hence, $(f_n(x))_{n\in\mathbb{N}}$ is a Cauchy sequence in the complete metric space $N$. Therefore, define the mapping $f:M\to N$ such that for all $x\in M\setminus Z$ it satisfies $f(x)\coloneqq \lim_{n\to\infty} f_n(x)$ and $f|_Z\equiv z_0$. Then, $f$ belongs to $\cLs(M,N)$, by \cref{prop:measurability_pointwise_limit}. Furthermore, we have for all $n \in \mathbb{N}^*$ and $x\in M\setminus Z$ that 
    \begin{align*}
    d_N(f(x),f_n(x)) &= \lim_{m\to \infty}
        d_N(f_m(x),f_n(x)) \leq \lim_{m\to\infty}\sum_{k=n}^{m-1} d_N(f_k(x), f_{k+1}(x))= \lvert g(x) - g_{n-1}(x)\rvert \leq \lvert g(x)\rvert.
    \end{align*}
    Hence, $f$ even belongs to $\cLph(M,N)$ and, by Lebesgue's dominated convergence theorem \cite[Theorem~2.4.5]{cohn2013measure}, $f_n\to f$ in $\cLph(M,N)$ as $n\to\infty$. 
\end{proof}

The reverse implication can be proved when restricting to a measure $\mu_M$ that is not purely infinite or nontrivial depending on the value of $p\in [1,\infty]$.

\begin{proposition}[Completeness of the target space]
\label{prop:completeness_target}
Let $h\in L^0(M,N)$. Then, the following assertions hold:
\begin{enumerate}[label=(\roman*)]
\item for every $p\in [1,\infty)$, if $\mu_M$ is not purely infinite and $\Lph(M,N)$ is complete, then $N$ is complete.
\item if $\mu_M$ is nontrivial and $L^\infty_h(M,N)$ is complete, then $N$ is complete.
\end{enumerate}
\end{proposition}

\begin{proof}
    In either case, \cref{cor:closed_embedding_target} shows that $N$ is isometric, up to a scaling factor in case (i), to a closed subset of the corresponding nonlinear Lebesgue space. Therefore, $N$ is complete as a closed subset of a complete metric space \cite[(3.14.5)]{dieudonne1960treatise}. 
\end{proof}
\begin{remark}[Why it fails when $\mu_M$ is purely infinite or trivial]
    When $\lvert N\vert =1$, the space $\Lph(M,N)$ is trivial (\cref{prop:trivial_lebesgue}), hence complete, but so is $N$, since it contains a single element. Thus, the conclusion of \cref{prop:completeness_target} remains true in this case, independently of the measure $\mu_M$. Now, suppose that $\lvert N\rvert > 1$. If $p\in [1,\infty)$ and $\mu_M$ is purely infinite, then $\Lph(M,N)$ is trivial by assertion (i) of \cref{prop:trivial_lebesgue}, and is therefore complete regardless of whether $N$ is complete. Similarly, if $p=\infty$ and $\mu_M$ is trivial, then $L^\infty_h(M,N)$ is trivial by assertion (ii) of \cref{prop:trivial_lebesgue}, and is therefore complete regardless of whether $N$ is complete. In both cases, $N$ may be chosen incomplete, showing that the corresponding assumptions on $\mu_M$ in \cref{prop:completeness_target} cannot be omitted.
\end{remark}

As a consequence of \cref{prop:completeness_lebesgue,prop:completeness_target}, the completeness of nonlinear Lebesgue spaces is entirely characterized by the completeness of their target space under the corresponding measure-theoretic assumptions.

\begin{theorem}[Characterization of completeness]
\label{th:charact_completeness_lebesgue}
Let $h\in L^0(M,N)$. Then, the following assertions hold:
\begin{enumerate}[label=(\roman*)]
    \item for every $p\in [1,\infty)$, if $\mu_M$ is not purely infinite, then $\Lph(M,N)$ is complete if and only if $N$ is complete.
    \item if $\mu_M$ is nontrivial, then $L^\infty_h(M,N)$ is complete if and only if $N$ is complete.
\end{enumerate}
\end{theorem}

\begin{remark}[Related results in the literature]
    A similar result is mentioned without proof by L.~Ambrosio, N.~Gigli and G.~Savaré \cite[Section~5.4]{ambrosio2005gradient} in the case of nonlinear Lebesgue spaces with constant base mapping, probability Borel measure $\mu_M$, and separable, in the metric sense, base and target spaces. 
\end{remark}

We now move on to the characterization of separability in nonlinear Lebesgue spaces.

\subsection{Characterization of separability}
\label{sec:separability}

Consider the notion of ($\mu_M$-essentially) countably generated measure spaces.
\begin{definition}[Countably generated measurable spaces]
\label{def:separable_measurable_space} A measure space $(M,\Sigma_M,\mu_M)$ is called:
\begin{enumerate}[label=(\roman*)]
    \item 
    \emph{countably generated} if there exists a countable subfamily $\mathcal{C}$ of $\Sigma_M$ such that $\Sigma_M = \sigma(\mathcal{C})$.
    \item \emph{$\mu_M$-essentially countably generated} if there exists a countable subfamily $\mathcal{C}$ of $\Sigma_M$ such that for all $A\in \Sigma_M$ there exists $A'\in \sigma(\mathcal{C})$ satisfying $\mu_M(A\Delta A') = 0$.
\end{enumerate}
\end{definition}
\begin{remark}[Related definitions in the literature]
    Unlike T.~Hytönen, J.~van Neerven, M.~Veraar and L. Weis \cite[Definition~1.2.27~(b)]{hytonen2016analysis}, we do not require the sets contained in the generating subfamily $\mathcal{C}$ to have finite $\mu_M$-measure in the definition of $\mu_M$-essentially countably generated measure spaces, hence we do not use the \enquote{$\mu_M$-countably generated} terminology to differentiate our definition from theirs.
\end{remark}

Then, nonlinear Lebesgue spaces inherit separability from their target and base spaces.

\begin{proposition}[Separability of nonlinear Lebesgue spaces]
\label{prop:separability_lebesgue_general}
    Let $h\in L^0(M,N)$ and $p\in [1,\infty)$. Suppose that $\Lph(M,N)$ is trivial (\cref{prop:trivial_lebesgue}) or that $N$ is separable and that we have a disjoint decomposition $M=M_0\cup M_1$ in $\Sigma_M$ such that $\Lph(M_0,N)$ is trivial on $(M_0,\Sigma_M|_{M_0},\mu_M|_{M_0})$ and $(M_1,\Sigma_M|_{M_1},\mu_M|_{M_1})$ is $\mu_M$-essentially 
    countably generated and $\sigma$-finite. 
    Then, $\Lph(M,N)$ is separable.
\end{proposition}

\begin{remark}[On the assumption of triviality]
    Recall that, for any choice of measure space $(M,\Sigma_M,\mu_M)$, we have, by \cref{prop:trivial_lebesgue}, that $\Lph(M,N)$ is trivial if and only if $\mu_M$ is purely infinite or $\lvert N\rvert = 1$.
\end{remark}

\begin{proof}[Proof of \cref{prop:separability_lebesgue_general}] The proof is partly inspired by the proofs of \cite[Theorem~4.13]{brezis2011functional} and \cite[Proposition~3.4.5]{cohn2013measure}, which both treat the real-valued case.

\noindent Let $h\in \cLs(M,N)$ and $p \in [1,\infty)$. When $\Lph(M,N)$ is trivial, that is, $\Lph(M,N)=\left\{[h]\right\}$, it is also separable as it contains a single element. Therefore, assume that $\Lph(M,N)$ is nontrivial, that $N$ is separable and that we have a disjoint decomposition $M=M_0\cup M_1$ in $\Sigma_M$ such that $\Lph(M_0,N)$ is trivial on $(M_0,\Sigma_M|_{M_0},\mu_M|_{M_0})$ and $(M_1,\Sigma_M|_{M_1},\mu_M|_{M_1})$ is $\mu_M$-essentially 
    countably generated and $\sigma$-finite. Now, let $f\in \cLph(M,N)$, $\varepsilon > 0$, $z_0\in N$ and $Q\coloneqq \left\{y_n\in N:n\in\mathbb{N}\right\}$ a countable dense subset of $N$.
In addition, since $\mu_M$ is $\sigma$-finite, there exists a collection $(M_n)_{n\in\mathbb{N}}$ of disjoint measurable sets such that $\mu_M(M_n)<\infty$ and $M=\cup_{n\in\mathbb{N}}M_n$. Also, define $H_n \coloneqq \left\{ x\in M : d_N(z_0, h(x))\leq n\right\}$, so that $M = \cup_{n\in\mathbb{N}} H_n$. 

\par\medskip \noindent \emph{Step 1 (Show that the restriction to $M_1$ is an isometry):} The fact that $\Lph(M_0,N)$ is trivial yields that $f|_{M_0}\sim  h|_{M_0}$. Therefore, the restriction mapping $f\in \cLph(M,N)\mapsto f|_{M_1}\in \cLph(M_1,N)$ is an isometry, so that we may assume $M=M_1$ in the rest of the proof.

\par\medskip \noindent \emph{Step 2 (Identification of the candidate countable dense set $\mathcal{Q}$):} Since $(M,\Sigma_M,\mu_M)$ is $\mu_M$-essentially countably generated, we can choose a countable subfamily $\mathcal{C}$ of $\Sigma_M$ containing both $(M_n)_{n\in\mathbb{N}}$ and $(H_n)_{n\in\mathbb{N}}$ such that for all $A\in \Sigma_M$ there is $A'\in \sigma(\mathcal{C})$ satisfying $\mu_M(A\Delta A') = 0$ and denote $\mathcal{A}$ the algebra generated by $\mathcal{C}$. Elements of $\mathcal{A}$ are generated by finite unions of elements of $\mathcal{C}$ and their complements, so $\mathcal{A}$ is countable (see the proof of \cite[Proposition~3.4.5]{cohn2013measure}). Now, define the set $\mathcal{Q}$ of mappings $f$ in $\cLph(M,N)$ such that there exists $n \in \mathbb{N}$ and a collection of pairs $((y_i,A_i))_{0\leq i\leq n}$ in $ (Q\times \mathcal{A})^{n+1}$ satisfying $A_i \cap A_j = \emptyset$ when $i\neq j$ and $f|_{A_i}\equiv y_i$ for all $0\leq i\leq n$ as well as $f|_{M\setminus A}= h|_{M\setminus A}$ with $A\coloneqq\cup_{i=0}^n A_i$. We then claim that $\mathcal{Q}$ is countable and dense in $\cLph(M,N)$. The fact that $\mathcal{Q}$ is countable follows from its construction which consists in choosing finitely many pairs of measurable sets in $\mathcal{A}$ and values in $Q$ to alter the base mapping $h$ on these measurable sets using the associated selected values. Thus, there is an injective mapping from $\mathcal{Q}$ to the countable set $\cup_{n\in \mathbb{N}} (Q\times \mathcal{A})^{n+1}$, hence $\mathcal{Q}$ is itself (at most) countable. 

\par\medskip \noindent \emph{Step 3 (Approximation by an almost simple mapping):} By \cref{prop:density_almost_simple} and \cref{rem:range_set_density_simple}, there exist $g\in \mathcal{E}_h(M,N)$, $B\in \mathcal{F}_{\mu_M}$ and $n_0\in \mathbb{N}$ such that $D_p(f,g) < \varepsilon/3$, $g|_B\in \mathcal{E}(M,N)$, $g|_{M\setminus B} \sim h|_{M\setminus B}$ and $g(B) = \left\{y_n: 0\leq n\leq n_0\right\}$. Also, define for all $0\leq n\leq n_0$ the measurable set $B_n \coloneqq g|_B^{-1}(\left\{y_n\right\})$, which, up to modifications of $g$ on a $\mu_M$-null set, can be assumed to belong to $\sigma(\mathcal{C})$. At this point, the $B_n$'s might not belong to $\mathcal{A}$, so that the key argument is now to approximate them by elements of $\mathcal{A}$ to produce an approximation mapping which belongs to $\mathcal{Q}$. The issue is that if the approximation sets go out of the $B_n$'s and the approximation mappings are set to constant values outside the $B_n$'s, the base mapping $h$, and hence $g$, might not be bounded outside the $B_n$'s, so that, without being careful in their construction, the approximation mappings could be at an infinite $\mathcal{L}^p$ distance from $g$. Thus, we need to restrict to subsets where $h$ is bounded, while keeping a good approximation.

\par\medskip \noindent \emph{Step 4 (Approximation by a mapping that differs from $h$ on a set where $h$ is bounded):} Then, define for all $n\in \mathbb{N}$ 
$$g_n(x)\coloneqq \begin{cases}
    g(x),&\text{if $d_N(z_0,h(x)) \leq n$}\\
    h(x),&\text{otherwise}
\end{cases}.$$
Then, since $(H_n)_{n\in\mathbb{N}}$ is an increasing sequence of sets in $M$, we have for all $x\in M$ that $g_n(x)\to g(x)$ as $n\to \infty$ and $d_N(g_n(x),g(x))^p \leq d_N(h(x),g(x))^p$. Hence, by Lebesgue's dominated convergence theorem \cite[Theorem~2.4.5]{cohn2013measure}, $D_p(g_n,g)^p\to 0$ as $n\to \infty$, so that we can pick $n_1\in \mathbb{N}$ such that $D_p(g_{n_1},g) < \varepsilon /3$. 

\par\medskip \noindent \emph{Step 5 (Approximation by an element of $\mathcal{Q}$):} Recall that 
$H_{n_1}= \left\{x\in M : d_N(z_0, h(x))\leq n_1\right\}\in \mathcal{A}$ 
and define 
$$C \coloneqq \max \left\{\max_{0\leq i<j\leq n_0} d_N(y_i,y_j), \max_{0\leq n\leq n_0}\{d_N(y_n,z_0) + n_1\}\right\}.$$
Also, define for all $0\leq n\leq n_0$ the set $B_n'\coloneqq B_n\cap H_{n_1}\subset B$, so that $\mu_M(B_n') < \infty$. Thus, together with the fact that $M=\cup_{n\in\mathbb{N}} M_n$ with $M_n\in \mathcal{A}$, we have, by \cite[Lemma~3.4.7]{cohn2013measure}, that there exists $A_n \in \mathcal{A}$ such that $\mu_M(B_n' \Delta A_n) < \varepsilon^p/(3\,(n_0+1)^{1/p}\,C)^p $, which can be made disjoint and contained in $H_{n_1}$. Also, define the mapping $\tilde{f}: M\to N$ such that for all $0\leq n\leq n_0$ it satisfies $\tilde{f}|_{A_n} \equiv y_n$ and $\tilde{f}|_{M\setminus A} = h|_{M\setminus A}$ with $A\coloneqq \cup_{n=0}^{n_0} A_n$. Then, by construction, $\tilde{f}$ belongs to $\mathcal{Q}$ and we have 
$$ D_p(\tilde{f}, g_{n_1})^p \leq C^p \sum_{n=0}^{n_0} \mu_M(B_n'\Delta A_n) < (\varepsilon/3)^p.$$
Finally, we get, by the triangle inequality, that
\begin{equation*}
    D_p(f,\tilde{f})\leq D_p(f,g) + D_p(g,g_{n_1}) + D_p(g_{n_1},\tilde{f})  < \varepsilon.\qedhere
\end{equation*}
\end{proof}
\begin{remark}[Related results in the literature]
    A similar result was proved by M.~Bauer, F.~Mémoli, T.~Needham and M.~Nishino \cite[Proposition~8]{bauer2025z} in the case where $h$ is constant, $(M,d_M)$ is a separable metric space (required in \cite[Proposition~8]{bauer2025z}), $\mu_M$ is a finite Borel measure (required in \cite[Remark~7]{bauer2025z}), and $N$ is separable. Since the Borel $\sigma$-algebra of a separable metric space is countably generated (see \cref{prop:borel_separable} for a proof), their result comes as a corollary of \cref{prop:separability_lebesgue_general}. 
\end{remark}

Then, the reverse implication can be proved by restricting to nonlinear Lebesgue spaces that are nontrivial.

\begin{proposition}[Separability of the target and base spaces]
\label{prop:separability_target_space}
    Let $h\in L^0(M,N)$ and $p\in [1,\infty)$. 
    Suppose that $\Lph(M,N)$ is nontrivial (\cref{prop:trivial_lebesgue}) and 
    separable. Then, the following assertions hold:
    \begin{enumerate}[label=(\roman*)]
        \item $N$ is separable. 
        \item we have a disjoint decomposition $M=M_0\cup M_1$ such that $\Lph(M_0,N)$ is trivial on $(M_0,\Sigma_M|_{M_0},\mu_M|_{M_0})$ and $(M_1,\Sigma_M|_{M_1},\mu_M|_{M_1})$ is $\mu_M$-essentially countably generated and $\sigma$-finite.
    \end{enumerate}
    In addition, $\Lph(M,N)\cong\Lph(M_1,N)$ isometrically.
\end{proposition}

\begin{remark}[Why it fails when $\Lph(M,N)$ is trivial]
    When $\Lph(M,N)$ is trivial, that is, $\Lph(M,N)=\left\{h\right\}$, it is separable regardless of the nature of the target and base spaces.
\end{remark}
\begin{remark}[On the triviality of $\Lph(M_0,N)$]
Since $\Lph(M,N)$ is assumed nontrivial, we already know that $\lvert N\rvert > 1$, so that, recalling \cref{prop:trivial_lebesgue}, $\Lph(M_0,N)$ being trivial becomes equivalent to $(M_0,\Sigma_M|_{M_0},\mu_M|_{M_0})$ being purely infinite.
\end{remark}
\begin{proof}[Proof of \cref{prop:separability_target_space}] The proof is partly inspired by the proof of the same implication in \cite[Proposition~1.2.29]{hytonen2016analysis}, which deals with the Banach space-valued case. However, the following proof does not rely on the density of simple mappings and the linear structure of the target space to show that the base space is $\mu_M$-essentially countably generated. 
    
    \noindent Let $h\in \cLs(M,N)$, $p\in [1,\infty)$ and $z_0\in N$.
    Since $\Lph(M,N)$ is nontrivial, we know, by \cref{prop:trivial_lebesgue}, that $\lvert N\rvert > 1$ and that $\mu_M$ is not purely infinite. 

    \begin{enumerate}[label=(\roman*),wide]
        \item 
        By \cref{cor:closed_embedding_target}, $N$ is isometric to a subset of the separable set $\Lph(M,N)$, hence is separable \cite[(3.10.9)]{dieudonne1960treatise}. 
    \item 
   Let $\mathcal{G}\coloneqq \varphi_h(\cLph(M,N))\subset \mathcal{L}^p(M,\mathbb{R}_+)$, which is separable as the continuous (\cref{prop:continuity_restriction}) image of a separable set. In particular, if $\mathcal{Q}$ is a countable dense subset of $\cLph(M,N)$, $\varphi_h(\mathcal{Q})$ is countable and dense in $\mathcal{G}$. Now, define for all $(f,n)\in \cLph(M,N)\times \mathbb{N}^*$ the measurable set $E_{f,n} \coloneqq \{x\in M: \varphi_h(f)(x) \geq n^{-1}\}$ and note that $\mu_M(E_{f,n})<\infty$, since $f\in \cLph(M,N)$. Then, define $\mathcal{A}\coloneqq \sigma(\{E_{f,n}: (f,n)\in \mathcal{Q}\times \mathbb{N}^*\})$, which is countably generated.
    
    \par\medskip \noindent \emph{Step 1 (Construct a $\mu_M$-essentially generated and $\sigma$-finite measure subspace):} Let $M_1 \coloneqq  \cup_{(f,n)\in\mathcal{Q}\times \mathbb{N}^*}E_{f,n}$. Then, $(M_1,\Sigma_1,\mu_1)$, where $\Sigma_1 \coloneqq \Sigma_M|_{M_1}$ and $\mu_1\coloneqq\mu_M|_{M_1}$, is $\sigma$-finite. We claim that it is also $\mu_M$-essentially countably generated. Indeed, let $A_1\in \Sigma_1$ and define for all $(f,n)\in \mathcal{Q}\times \mathbb{N}^*$ the set $A_{f,n} \coloneqq A_1\cap E_{f,n}$. 
    Then, define the mapping $g: M\to N$ as
    $g|_{A_{f,n}} = f$ 
    and $g|_{M\setminus A_{f,n}} = h$, which belongs to $\cLph(M,N)$. Thus, $\varphi_h(g)$ belongs to $\mathcal{G}$ and vanishes outside $A_{f,n}$. By density of $\varphi_h(\mathcal{Q})$ in $\mathcal{G}$, there exists a sequence $(g_k)_{k\in\mathbb{N}}$ in $\mathcal{Q}$ such that $\lVert \varphi_h(g) - \varphi_h(g_k)\rVert_{p,\mu_M}^p \to 0 $ as $k\to \infty$. Up to the extraction of a subsequence (see the proof of \cite[Theorem~4.8]{brezis2011functional}), we can assume that for all $k\in\mathbb{N}$ it holds $\lVert \varphi_h(g) -\varphi_h(g_k)\rVert_{p,\mu_M}^p \leq 2^{-k}$. 
    Note that, for all $k\in\mathbb{N}$ and $x\in A_{f,n}\setminus E_{g_k,n+1}$ the following holds
    $$\lvert \varphi_h(g)(x) - \varphi_h(g_k)(x)\rvert  \geq  n^{-1} - (n+1)^{-1}.$$
    Hence, we get that for all $k\in \mathbb{N}$ that
    \begin{align*}
    (n(n+1))^{-p} \mu_1(A_{f,n}\setminus E_{g_k,n+1}) &\leq\int_M \lvert \varphi_h(g)(x) - \varphi_h(g_k)(x)\rvert^p\dif\mu_M(x)\leq 2^{-k},
    \end{align*}
    so that 
    $\lim_{k\to \infty} \mu_1(A_{f,n}\setminus E_{g_k,n+1})= 0$. Therefore, defining $B_{f,n} \coloneqq \cap_{l\in\mathbb{N}}\cup_{k\geq l} E_{g_k,n+1}\in \mathcal{A}$ and, using the continuity from below of $\mu_M$ \cite[Proposition~1.2.5~(a)]{cohn2013measure}, we get that
    \begin{align*}\mu_1(A_{f,n} \setminus B_{f,n}) &= \mu_1(\cup_{l\in \mathbb{N}} \cap_{k\geq l} (A_{f,n}\setminus E_{g_k,n+1})) = \lim_{l\to \infty} \mu_1(\cap_{k\geq l} (A_{f,n}\setminus E_{g_k,n+1}))\leq \lim_{l\to \infty} \mu_1(A_{f,n}\setminus E_{g_l,n+1}) = 0.
    \end{align*}
    Furthermore, using Markov's inequality \cite[Proposition~2.3.10]{cohn2013measure}, we also get for all $k\in \mathbb{N}$ that
    \begin{align*}
    (n+1)^{-p} \mu_1(E_{g_k,n+1}\setminus A_{f,n})&\leq \int_{M\setminus A_{f,n}} \lvert \varphi_h(g_k)(x)\rvert ^p\dif\mu_M(x)\leq \int_M \lvert \varphi_h(g)(x) - \varphi_h(g_k)(x)\rvert^p\dif\mu_M(x)\leq 2^{-k},
    \end{align*}
    which yields that
    $\sum_{k\in\mathbb{N}} \mu_1(E_{g_k,n+1}\setminus A_{f,n})< \infty$ and, by Borel-Cantelli's lemma \cite[Proposition~10.2.2~(a)]{cohn2013measure}, $\mu_1(\cap_{l\in\mathbb{N}}\cup_{k\geq l} (E_{g_k,n+1}\setminus A_{f,n})) =0$. 
    Therefore, we have 
    \begin{align*}\mu_1(B_{f,n}\setminus A_{f,n}) = \mu_1(\cap_{l\in\mathbb{N}} \cup_{k\geq l} (E_{g_k,n+1}\setminus A_{f,n})) = 0.
    \end{align*}
    Thus, $\mu_1(A_{f,n}\Delta B_{f,n}) = 0$ and, since we can find such a set for all $(f,n)\in \mathcal{Q}\times \mathbb{N}^*$, we can define $B_1\coloneqq \cup_{(f,n)\in\mathcal{Q}\times \mathbb{N}^*} B_{f,n}$ which belongs to $\mathcal{A}$ and is such that \begin{align*}\mu_1(A_1\Delta B_1)&\leq \mu_1(\cup_{(f,n)\in \mathcal{Q}\times \mathbb{N}^*} (A_{f,n}\Delta B_{f,n}))\leq \sum_{(f,n)\in\mathcal{Q}\times \mathbb{N}^*} \mu_1(A_{f,n}\Delta B_{f,n})=0,
    \end{align*}
    which proves that $(M_1,\Sigma_1,\mu_1)$ is $\mu_M$-essentially countably generated.
    
    \par\medskip \noindent \emph{Step 2 (Show that the nonlinear Lebesgue space defined on the complement measure subspace is trivial):} Now, define $M_0 \coloneqq M\setminus M_1$. Then, we claim that $\Lph(M_0,N)$ is trivial on $(M_0,\Sigma_0,\mu_0)$, where $\Sigma_0 \coloneqq \Sigma_M|_{M_0}$ and $\mu_0\coloneqq\mu_M|_{M_0}$. Thus, let $f\in \cLph(M_0,N)$. By \cref{prop:lebesgue_differ_from_base_bounded}, $f$ differs from $h$ on a $\sigma$-finite subset $A_0\in \Sigma_0$. 
Precisely, we have $A_n \coloneqq \{x\in M_0: d_N(f(x),h(x)) \geq n^{-1} \} \subset A_0$ for every $n\in\mathbb{N}^*$, so that $A_0 = \cup_{n\in\mathbb{N}^*} A_n$, and define the mapping $\tilde{g}_n: M\to N$ such that $\tilde{g}_n|_{A_n} = f|_{A_n}$ and $\tilde{g}_n|_{M\setminus A_n} = h|_{M\setminus A_n}$. Thus, $\varphi_h(\tilde{g}_n)$ belongs to $\mathcal{G}$ and vanishes outside $A_n$. By density of $\varphi_h(\mathcal{Q})$ in $\mathcal{G}$, there exists a sequence $(g_k)_{k\in\mathbb{N}}$ in $\mathcal{Q}$ such that $\lVert \varphi_h(\tilde{g}_n) - \varphi_h(g_k)\rVert_{p,\mu_M}^p \to 0 $ as $k\to \infty$. Since the $\varphi_h(g_k)$'s vanish outside $M_1$, the latter yields $\mu_0(A_n) = \mu_M(A_n\cap M_0) = 0$. Finally, $A_0= \cup_{n\in\mathbb{N}^*} A_n$, thus, using the countable subadditivity of $\mu_0$, we get $\mu_0(A_0) \leq \sum_{n\in\mathbb{N}^*} \mu_0(A_n) = 0$, that is, $f$ differs from $h$ on a $\mu_0$-null set on $M_0$. Therefore, $\Lph(M_0,N)$ is trivial. 

\par\medskip \noindent \emph{Step 3 (Show that the restriction to $M_1$ is an isometry):} The restriction mapping $f \in \Lph(M,N)\mapsto f|_{M_1}\in \Lph(M_1,N)$ is an isometry since $f|_{M_0}\sim h|_{M_0}$ for all $f\in \Lph(M,N)$.
\qedhere
    \end{enumerate}
\end{proof}

\begin{remark}[Related results in the literature]
    A similar result was proved by T.~Hytönen, J.~van Neerven, M.~Veraar and L.~Weis \cite[Proposition~1.2.29]{hytonen2016analysis} in the linear case where $N$ is additionally a Banach space.
\end{remark}

\cref{prop:separability_lebesgue_general,prop:separability_target_space} then yield a characterization of separability in nonlinear Lebesgue spaces.

\begin{theorem}[Characterization of separability]
\label{th:charact_separability_lebesgue_general}
    Let $h\in L^0(M,N)$ and $p\in [1,\infty)$. Suppose that $\Lph(M,N)$ is nontrivial (\cref{prop:trivial_lebesgue}). 
    Then, the following assertions are equivalent:
    \begin{enumerate}[label=(\roman*)]
        \item $\Lph(M,N)$ is separable.
        \item $N$ is separable and we have a disjoint decomposition $M=M_0\cup M_1$ in $\Sigma_M$ such that $\Lph(M_0,N)$ is trivial on $(M_0,\Sigma_M|_{M_0},\mu_M|_{M_0})$ and $(M_1,\Sigma_M|_{M_1},\mu_M|_{M_1})$ is $\mu_M$-essentially countably generated and $\sigma$-finite. 
    \end{enumerate}
    If these equivalent conditions hold, we have $\Lph(M,N)\cong\Lph(M_1,N)$ isometrically.
\end{theorem}

\begin{remark}[Related results in the literature]
    This characterization of separability is well-known in the context of linear Lebesgue spaces with constant base mapping (see \cite[Proposition~1.2.29]{hytonen2016analysis}, \cite[Exercise~4.7.63]{bogachev2007measure} or \cite[Exercise~365X~(p)]{fremlin2000measure}).
\end{remark}

\section{Density of simple and countably valued mappings}
\label{sec:density_simple}

As in the linear case, the space of simple mappings is dense in nonlinear Lebesgue spaces under some conditions on the base mapping and/or the measure.

\subsection{\texorpdfstring{Case $p\in [1,\infty)$}{Case p in [1,infty)}}

First, this holds for $p\in [1,\infty)$ with conditions on the base mapping and/or on $\mu_M$.

\begin{theorem}[Density of simple mappings for $p\geq 1$]
\label{th:density_simple}
Let $p\in [1,\infty)$. Suppose that one of the following conditions holds:
\begin{enumerate}[label=(\roman*)]
    \item \label{cond:finite} $h\in \Lb(M,N)$ and $\mu_M$ is finite.
    \item \label{cond:simple} $h\in E(M,N)$.
\end{enumerate}
Then, $E(M,N)\cap \Lph(M,N)$ is a dense subspace of $\Lph(M,N)$.
\end{theorem}
\begin{remark}[Explicit expression of the dense subspace]
    $E(M,N)\cap \Lph(M,N) = E_h(M,N)$ (\cref{def:almost_simple}), when $h\in E(M,N)$, (see \cref{prop:simple_lebesgue_are_alsmot_simple} for a proof) and $E(M,N)\cap \Lph(M,N) = E(M,N)$, when $h\in \Lb(M,N)$ and $\mu_M$ is finite.
\end{remark}
\begin{proof}[Proof of \cref{th:density_simple}]
    Let $p\in [1,\infty)$.
    \begin{enumerate}[label=(\roman*),wide]
        \item Since $\mu_M$ is finite and $h$ is bounded, we have, by \cref{prop:indep_def}, that $\Lph(M,N) = L^p_{h'}(M,N)$ with $h'\equiv z_0\in N$. Then, by \cref{prop:density_almost_simple}, $E_{h'}(M,N)\subset E(M,N)$ is a dense subspace of $L^p_{h'}(M,N)$, hence $E(M,N)$ is a dense subspace of $\Lph(M,N)$.
        \item This is a consequence of \cref{prop:density_almost_simple} as, when $h\in E(M,N)$, $E_h(M,N)$ is a subset of $E(M,N)\cap \Lph(M,N)$.\qedhere
    \end{enumerate}
\end{proof}
\begin{remark}[Related results in the literature]
\label{rem:cor_neumayer}
    A similar result was proved by S.~Neumayer, J.~Persch and G.~Steidl \cite[Lemma~2.1]{neumayer2018} in the case where $h$ is constant, $M$ is an open, bounded, connected (Lipschitz) domain of $\mathbb{R}^n$, $n\geq 2$, $\mu_M=\mathscr{L}^n$ and $N$ is a locally compact Hadamard space (see \cite[Definition~1.2.3]{bacak2014convex} for the definition of Hadamard spaces). In particular, a Hadamard space is, by definition, a geodesic space, hence a length space (see \cite{burago2022course} for the definition of length and geodesic spaces). Since a locally compact length space is separable (see \cite[Appendix A,~LEMMA,~p.~460]{spivak1979comprehensive} for the proof that any connected, locally compact and paracompact space is the countable union of compact sets, hence is separable), their result comes as a corollary of \cref{th:density_simple}. Their proof relies on the fact that, under these assumptions, $N$ is boundedly compact (see \cite[Theorem~2.5.28]{burago2022course} for a proof) to produce a \enquote{sufficiently large} compact, hence separable, subset of $N$, which we avoid by exploiting the fact that involved mappings have separable range. Unlike \cref{th:density_simple}, the local compactness assumption enforces a finite-dimensional constraint when $N$ is a Banach space \cite[Chapter II,  \S 3, Corollary 3.15]{lang2012real}, so that it is not a desirable assumption in our general setting.
\end{remark}

We now discuss the sharpness of the assumptions used in \cref{th:density_simple} by providing counterexamples. The reader may skip to \cref{sec:separabiliy_infty} in the first reading.
First, we begin by discussing the condition \ref{cond:finite} of \cref{th:density_simple} and start by giving a counterexample when the base mapping is unbounded.

\begin{cexample}[When $\mu_M$ is finite but $h$ is unbounded]
Let $p\in [1,\infty)$. Suppose that $M=(0,1]$ is equipped with the $\sigma$-algebra $\Sigma_M= \mathcal{B}((0,1])$ and the finite measure $\mu_M=\mathscr{L}^1|_{(0,1]}$, and that $N=\mathbb{R}_+$ is equipped with its standard metric. Then, let $h: x\mapsto x^{-\frac{1}{p}}$, which is bounded on any measurable subset of $(0,1]$ that does not contain $0$ in its closure, but is unbounded on $(0,1]$. Note that $f\in \cLph((0,1],\mathbb{R}_+)$ if and only if the mapping $x \mapsto x^{-\frac{1}{p}} - f(x)$ belongs to $\mathcal{L}^p((0,1],\mathbb{R})$. Hence, no simple mapping belongs to $\cLph((0,1],\mathbb{R}_+)$ since such mappings are bounded. Indeed, if $g\in \mathcal{E}((0,1],\mathbb{R}_+)$, its boundedness implies that $g(x)/h(x) = x^{1/p}g(x)\to 0$ as $x\to 0$, so that $\lvert g(x) - h(x)\rvert^p \sim h(x)^p=x^{-1}$ as $x\to 0$. Therefore, we get that
$$D_p(g,h)^p = \int_0^1 \lvert g(x) - h(x)\vert^p\dif x = \infty.$$ 
\end{cexample}

A counterexample can also be provided when the condition on $\mu_M$ is relaxed to $\sigma$-finiteness.

\begin{cexample}[When $h$ is bounded but $\mu_M$ is $\sigma$-finite]
    Let $p\in [1,\infty)$. Suppose that $M=\mathbb{R}_+$ is equipped with the $\sigma$-algebra $\Sigma_M =\mathcal{B}(\mathbb{R}_+)$ and the $\sigma$-finite measure $\mu_M=\mathscr{L}^1|_{\mathbb{R}_+}$, and that $N=\mathbb{R}$ is equipped with its standard metric. Then, let $h:x\mapsto (1+x)^{-1/p}$, which is bounded on $\mathbb{R}_+$ but not simple. In that case, no simple mapping belongs to $\cLph(\mathbb{R}_+,\mathbb{R})$. Indeed, suppose that $\mathcal{E}(\mathbb{R}_+,\mathbb{R})\cap \cLph(\mathbb{R}_+,\mathbb{R})\neq\emptyset$ and let $g\in \mathcal{E}(\mathbb{R}_+,\mathbb{R})\cap \cLph(\mathbb{R}_+,\mathbb{R})$, so that we can assume there exists a finite indexing set $I$ such that $g(\mathbb{R}_+)=\{y_i \in \mathbb{R}: i \in I\}$ and define $M_i\coloneqq g^{-1}(\{y_i\})$. For every $i\in I$ such that $y_i\neq 0$, choose $R_i > 0$ such that $h(x)\leq \lvert y_i\vert/2$ for all $x\geq R_i$. Then, $\lvert y_i - h(x)\vert\geq \lvert y_i \rvert/2$ for all $x\in M_i\cap [R_i,\infty)$. Therefore, $D_p(g,h)^p\geq (\lvert y_i\rvert/2)^p\mathscr{L}^1(M_i\cap[R_i,\infty))$. Hence, $\mathscr{L}^1(M_i) < \infty$ whenever $y_i\neq 0$ as $g\in\mathcal{L}^p_h(\mathbb{R}_+,\mathbb{R})$. Since $g$ is simple, the set $A\coloneqq \{x: g(x)\neq 0\}$ satisfies $A=\cup_{\{i\in I: y_i\neq 0\}} M_i$ and has finite Lebesgue measure. In addition, $h(x)^p = (1+x)^{-1} \leq 1$ for every $x\in \mathbb{R}_+$, so $\int_A h(x)^p\dif x\leq \mathscr{L}^1(A) < \infty$. However, $\int_{\mathbb{R}_+} h(x)^p\dif x = \int_0^\infty(1+x)^{-1}\dif x  =\infty$. Thus, $\int_{\mathbb{R}_+\setminus A} h(x)^p\dif x= \infty$ and, since $g|_{\mathbb{R}_+\setminus A}\equiv 0$, we get that  
    $$D_p(g,h)^p\geq \int_{\mathbb{R}_+\setminus A} \lvert g(x) - h(x)\rvert^p \dif x = \int_{\mathbb{R}_+\setminus A}h(x)^p\dif x = \infty,$$
    which yields a contradiction.
\end{cexample}

We now move on to the discussion on the condition \ref{cond:simple} of \cref{th:density_simple} by providing a counterexample when the base mapping is not simple but countably valued.

\begin{cexample}[When $\mu_M$ is $\sigma$-finite but $h$ is countably valued]
    Let $p\in [1,\infty)$. Suppose $M=\mathbb{R}_+$ is equipped with the $\sigma$-algebra $\Sigma_M = \mathcal{B}(\mathbb{R}_+)$ and the $\sigma$-finite measure $\mu_M = \mathscr{L}^1|_{\mathbb{R}_+}$, and that $N= \mathbb{R}$ is equipped with its standard metric. Then, define the sequence $(a_n)_{n\in\mathbb{N}}$ by $a_0 \coloneqq 0$ and $a_n\coloneqq \sum_{k=1}^n k$ for $n\geq 1$, so that the sequence $(I_n)_{n\in\mathbb{N}}$ defined as $I_n \coloneqq [a_n,a_{n+1})$ forms a partition of $\mathbb{R}_+$. Now, let $h: \mathbb{R}_+ \to \mathbb{R}$ be defined as $h|_{I_n} \coloneqq (n+1)^{-1/p}$ for all $n\in \mathbb{N}$, which is simple on every bounded measurable subset of $\mathbb{R}_+$, but is countably valued on $\mathbb{R}_+$. In that case, no simple mapping belongs to $\cLph(\mathbb{R}_+,\mathbb{R})$. Indeed, suppose that $\mathcal{E}(\mathbb{R}_+,\mathbb{R})\cap \cLph(\mathbb{R}_+,\mathbb{R})\neq \emptyset$ and let $g\in \mathcal{E}(\mathbb{R}_+,\mathbb{R})\cap \cLph(\mathbb{R}_+,\mathbb{R})$, so that we can assume there exists a finite indexing set $I$ such that $g(\mathbb{R}_+) = \left\{y_i \in \mathbb{R}:i\in I\right\}$ and define $M_i \coloneqq g^{-1}(\left\{y_i\right\})$. Then,
    $$D_p(g,h)^p = \sum_{i\in I} \sum_{n\in \mathbb{N}}  \lvert y_i - (n+1)^{-1/p}\rvert^p \mathscr{L}^1(M_i\cap I_n) < \infty.$$
    Since $\mathbb{R}_+ = \cup_{i\in I} M_i$, there exists $i_0\in I$ such that both $y_{i_0}\neq 0$ and $\mathscr{L}^1(M_{i_0}) = \infty$ otherwise $\mathscr{L}^1(\mathbb{R}_+)= \sum_{i\in I} \mathscr{L}^1(M_i)$ would be finite as a finite sum of finite values or $D_p(g,h) = \infty$ if we had $y_i \neq 0$ only for $M_i$'s of finite Lebesgue measure. In particular, this yields that the sequence of nonnegative reals $(\mathscr{L}^1(M_{i_0}\cap I_n))_{n\in\mathbb{N}}$ has a diverging sum.
    However, $\lvert y_{i_0} - (n+1)^{-1/p}\rvert^p \to \lvert y_{i_0}\rvert^p$ as $n\to \infty$, so that there exists $n_0\in \mathbb{N}$ such that for all $n\geq n_0$ we have the following inequality
    $$\lvert y_{i_0} - (n+1)^{-1/p}\rvert^p \mathscr{L}^1(M_{i_0}\cap I_n) \geq \lvert y_{i_0}\rvert^p/2\,\mathscr{L}^1(M_{i_0}\cap I_n).$$
    Hence, $\sum_{n\in \mathbb{N}} \lvert y_{i_0} - (n+1)^{-1/p}\rvert^p \mathscr{L}^1(M_{i_0}\cap I_n) = \infty$ and thus $D_p(g,h) = \infty$, which is absurd.
\end{cexample}

When relaxing the assumptions on the base mapping and/or the measure $\mu_M$, one can still retrieve the density of countably-valued mappings.

\begin{proposition}[Density of countably valued mappings for $p\geq 1$]
\label{prop:density_countably_valued_mappings}
    Let $h\in L^0(M,N)$ and $p\in [1,\infty)$. Suppose that one of the following conditions holds:
    \begin{enumerate}[label=(\roman*)]
        \item $\mu_M$ is $\sigma$-finite.
        \item $h$ is countably-valued.
    \end{enumerate}
    Then, the set of countably valued mappings that belong to $\cLph(M,N)$ is a dense subspace of $\cLph(M,N)$ and so determines a dense subspace of $\Lph(M,N)$.
\end{proposition}

\begin{remark}[On the necessity of choosing a $\sigma$-finite measure]
    We know from \cref{prop:lebesgue_differ_from_base_bounded} that Lebesgue mappings only differ from the base mapping on a $\sigma$-finite set $B$. However, this statement is not useful here since, when $h$ is not countably-valued, we cannot find a countably-valued mapping at a finite $\mathcal{L}^p$ distance from $h$ outside $B$ since $\mu_M(M\setminus B)$ could have no decomposition in measurable sets of finite $\mu_M$-measure. So, either $h$ should be assumed countably-valued itself or $\mu_M$ should be assumed $\sigma$-finite to ensure that $M$ could be decomposed in measurable sets of finite $\mu_M$-measure, so that a countably-valued approximation of $h$ would be at a finite $\mathcal{L}^p$ distance of $h$ on each of these measurable sets.
\end{remark}

\begin{proof}[Proof of \cref{prop:density_countably_valued_mappings}]
    Let $p\in [1,\infty)$ and $h\in \cLs(M,N)$.

    \begin{enumerate}[label=(\roman*),wide]
        \item Let $f\in \cLph(M,N)$ and $\varepsilon > 0$. Since $\mu_M$ is $\sigma$-finite, there exists a sequence $(M_n)_{n\in\mathbb{N}}$ of disjoint measurable sets with finite $\mu_M$-measure such that $M= \cup_{n\in\mathbb{N}} M_n$.
        By \cref{prop:countably_finite_radius}, we can find for all $n\in \mathbb{N}$ a countably valued measurable mapping $g_n$ satisfying $D_\infty(f,g_n) < \varepsilon/(2^{n+1}(1+\mu_M(M_n)))^{1/p}$.
        Now, define $g: M\to N $ such that $g|_{M_n} \coloneqq g_n|_{M_n}$ for all $n\in \mathbb{N}$, which is countably-valued since $g(M) = \cup_{n\in\mathbb{N}} g_n(M_n)$ is the countable union of countable sets. Hence, 
        \begin{align*}D_p(f,g)^p &= \sum_{n\in\mathbb{N}} \int_{M_n} d_N(f(x),g_n(x))^p\dif\mu_M(x)\leq \sum_{n\in\mathbb{N}} D_\infty(f,g_n)^p\mu_M(M_n) <  \varepsilon^p \sum_{n\in\mathbb{N}} 2^{-(n+1)}\frac{\mu_M(M_n)}{1+\mu_M(M_n)}\leq\varepsilon^p.
        \end{align*}
    
        \item This is a consequence of \cref{prop:density_almost_simple} as, when $h$ is countably-valued, so is any element of $E_h(M,N)$.\qedhere
    \end{enumerate}
\end{proof}

\subsection{\texorpdfstring{Case $p=\infty$}{Case p=infty}}
\label{sec:separabiliy_infty}

The density of simple mappings can also be proven in the case where $p=\infty$ but should be treated differently as it requires the stronger assumption that $N$ is boundedly compact, that is, closed and bounded sets are compact, yet weaker assumptions on $h$ and $\mu_M$.
\begin{proposition}[Density of simple mappings for $p=\infty$]
\label{prop:density_simple_infty}
    Let $h\in \Lb(M,N)$. Suppose that $N$ is boundedly compact. Then, $E(M,N)$ is a dense subspace of $L^\infty_h(M,N)$.
\end{proposition}
\begin{proof}
    This follows from a generalization of \cite[Proposition 3.4.2]{cohn2013measure} in the case $p=\infty$.

    \noindent Let $h\in \cLbs(M,N)$, $f\in \mathcal{L}^\infty_h(M,N)$, $z_0\in N$ and $\varepsilon > 0$. Then, by definition of the $\mu_M$-essential supremum, there exists $Z\in \Sigma_M$ such that $\mu_M(Z) = 0$ and for all $x\in M\setminus Z$ it holds that $d_N(f(x),h(x))\leq D_\infty(f,h) < \infty$. Thus, since $h$ is bounded, $f(M\setminus Z)$ is bounded and can be covered by a closed metric ball $C$, which is compact since $N$ is boundedly compact. Hence, we can pick $n_0\in \mathbb{N}$ and a finite collection $(y_i)_{0\leq i\leq n_0}$ of elements of $N$ such that the open metric balls centered on $y_i$ of radius $\varepsilon$, denoted $B_i$, cover the compact $C$. Then, the finite collection $(P_i)_{0\leq i\leq n_0}$ of disjoint Borel sets defined as $P_0\coloneqq B_0$ and for all $1\leq i\leq n_0$ as $P_i \coloneqq B_i\setminus \cup_{j = 0}^{i-1} B_j$ forms a partition of $C$. Define $\chi : C \to \{0,\ldots,n_0\}$ such that for all $y\in C$ it satisfies $y \in P_{\chi(y)}$ and define the mapping $g : M\to N$ as
    $$g(x) \coloneqq \begin{cases}y_{\chi \circ f(x)},&\text{if}\ f(x)\in C\\
    z_0,&\text{otherwise}
    \end{cases}.$$\\
    By construction, $g\in \mathcal{E}(M,N)$ and $D_\infty(f,g) <\varepsilon$.
\end{proof}

\begin{remark}[Related results in the literature]
    This is a well-known result in the linear case \cite[Proposition~3.4.2]{cohn2013measure}.
\end{remark}

As done in the case of $p\in [1,\infty)$, we now discuss the sharpness of the assumptions used in \cref{prop:density_simple_infty} by providing counterexamples. The reader may skip to \cref{sec:density_continuous} in the first reading.
First, we provide a counterexample when $N$ is not boundedly compact.

\begin{cexample}[When $N$ is not boundedly compact]
    Suppose $M=\mathbb{R}_+$ is equipped with the $\sigma$-algebra $\Sigma_M = \mathcal{B}(\mathbb{R}_+)$ and the measure $\mu_M = \mathscr{L}^1|_{\mathbb{R}_+}$, and that $N= \mathbb{H}$ is an infinite-dimensional separable Hilbert space, which is therefore not boundedly compact \cite[Theorem 6.5]{brezis2011functional}. Then, let $h\equiv 0_{\mathbb{H}}$, which is bounded. As a separable Hilbert space, $\mathbb{H}$ has an orthonormal basis $(e_n)_{n\in\mathbb{N}}$ \cite[Theorem~5.11]{brezis2011functional}. Then, define the sequence $(a_n)_{n\in\mathbb{N}}$ such that $a_0 \coloneqq 0$ and $a_n\coloneqq \sum_{k=1}^n k$ for $n\geq 1$, so that the sequence $(I_n)_{n\in\mathbb{N}}$ defined as $I_n \coloneqq [a_n,a_{n+1})$ forms a partition of $\mathbb{R}_+$. Then, let $f: \mathbb{R}_+ \to \mathbb{H}$ be defined as $f|_{I_n} \coloneqq e_n$ for all $n\in \mathbb{N}$ and note that $f\in \mathcal{L}^\infty_h(\mathbb{R}_+,\mathbb{H})$, since $\lVert f(x)\rVert_{\mathbb{H}}=1$ for every $x\in \mathbb{R}_+$. Let $g\in \mathcal{E}(\mathbb{R}_+,\mathbb{H})\cap \mathcal{L}_h^\infty(\mathbb{R}_+,\mathbb{H})$, so that we can assume there exists a finite indexing set $I$ such that $g(\mathbb{R}_+) = \left\{y_i\in \mathbb{H}: i\in I\right\}$ and define $M_i \coloneqq g^{-1}(\left\{y_i\right\})$. In addition, if $y_i^{(n)} \coloneqq \langle y_i, e_n\rangle_{\mathbb{H}}$, the Bessel-Parseval's identity yields that $(y_i^{(n)})_{n\in\mathbb{N}}$ belongs to $\mathcal{L}^2(\mathbb{N},\mathbb{R})$ \cite[Chapter 5, Remark 10, p.~144]{brezis2011functional}, hence $y_i^{(n)}\to 0$ as $n\to \infty$. Now, since $I$ is finite, there exists $i_0\in I$ such that $\mathscr{L}^1(M_{i_0}) = \infty$. Since each $I_n$ has finite Lebesgue measure, the set $J\coloneqq \{n\in\mathbb{N}: \mathscr{L}^1(M_{i_0}\cap I_n) > 0\}$ is infinite. Since $y_{i_0}^{(n)}\to 0$ as $n\to \infty$, there exists $n_0\in J$ such that $\lvert y_{i_0}^{(n_0)}\rvert < 1/2$. On the set $M_{i_0}\cap I_{n_0}$, which has positive Lebesgue measure, we have $g|_{M_{i_0}\cap I_{n_0}}\equiv y_{i_0}$ and $f|_{M_{i_0}\cap I_{n_0}}\equiv e_{n_0}$. Therefore, $$D_\infty(g,f)\geq \lVert y_{i_0} - e_{n_0}\rVert_{\mathbb{H}}\geq \lvert \langle y_{i_0} - e_{n_0},e_{n_0}\rangle_{\mathbb{H}}\rvert = \lvert y_{i_0}^{(n_0)} - 1\rvert > 1/2.$$
    Since this holds for every $g\in \mathcal{E}(\mathbb{R}_+,\mathbb{H})\cap \mathcal{L}^\infty_h(\mathbb{R}_+,\mathbb{H})$, no simple mapping can approximate $f$ in $\mathcal{L}^\infty_h(\mathbb{R}_+,\mathbb{H})$.
\end{cexample}

Now, we provide a counterexample when the base mapping is unbounded.

\begin{cexample}[When $h$ is unbounded]
    Suppose $M=[1,\infty)$ is equipped with the $\sigma$-algebra $\Sigma_M =\mathcal{B}([1,\infty))$ and the measure $\mu_M=\mathscr{L}^1|_{[1,\infty)}$, and that $N=\mathbb{R}$ is equipped with its standard metric, for which $\mathbb{R}$ is boundedly compact. Then, let $h:x\mapsto e^x$, which is bounded on every bounded subset of $[1,\infty)$. If $f\in \mathcal{L}^\infty_h([1,\infty),\mathbb{R})$, there exists $Z\in \mathcal{Z}_{\mu_M}$ such that $\lvert f(x) - h(x)\rvert \leq \lVert f - h\rVert_\infty < \infty$ for all $x\in [1,\infty)\setminus Z$. Therefore, $f(x)\geq e^x  - \lVert f - h\rVert_\infty$ for every $x\in [1,\infty)\setminus Z$. Since the set $[1,\infty)\setminus Z$ is unbounded, $f$ is unbounded and hence cannot be simple. Therefore, $\mathcal{E}([1,\infty),\mathbb{R})\cap \mathcal{L}^\infty_h([1,\infty),\mathbb{R}) = \emptyset$.
\end{cexample}

When relaxing the assumptions on the base mapping and $N$, one can still retrieve the density of the set of countably valued mappings in nonlinear $\infty$-Lebesgue spaces as an immediate consequence (by restriction) of \cref{prop:countably_finite_radius}.

\begin{proposition}[Density of countably valued mappings for $p=\infty$]
\label{cor:density_countably_valued_mappings}
    Let $h\in L^0(M,N)$. Then, the set of countably valued mappings that belong to $\mathcal{L}^\infty_h(M,N)$ is a dense subspace of $\mathcal{L}^\infty_h(M,N)$ and so determines a dense subspace of $L^\infty_h(M,N)$.
\end{proposition}

\begin{remark}[Boundedness of $\mathcal{L}^\infty$ mappings when $h$ is bounded]
    If $h$ is bounded, so are the elements of $\mathcal{L}_h^\infty(M,N)$ and thus the dense mappings (see \cite[Lemma~2.1.4]{hytonen2016analysis} for the linear case).
\end{remark}

\section{Density of continuous mappings}
\label{sec:density_continuous}

Another standard density result in the linear case is the density of compactly supported continuous mappings (see \cite[Proposition~7.4.3]{cohn2013measure}, \cite[Corollary~4.23]{brezis2011functional} or \cite[Lemma~1.2.31]{hytonen2016analysis}). The notion of compactly supported mapping does not make sense in our fully general setting, since $N$ is not necessarily a vector space. A slightly more general notion that fits this general setting is that of continuous mappings that take a constant value outside a compact set.

\begin{definition}[Continuous mappings]
\label{def:continuous_mappings}
    When $M$ is a topological space, define:
    \begin{enumerate}[label=(\roman*)]
        \item $\mathcal{C}(M,N)$ the space of \emph{continuous mappings}.
        \item $\Ccr(M,N)$ the space of \emph{continuous mappings with relatively compact range}, that is, all mappings $f\in \mathcal{C}(M,N)$ such that the closure of $f(M)$ in $N$ is compact.
        \item $\Cc(M,N)$ the space of \emph{continuous mappings that are constant outside a compact set}, that is, all mappings $f\in \mathcal{C}(M,N)$ such that there exist a compact $K\subset M$ and a point $z_0\in N$ for which $f(x) = z_0$ for all $x\in M\setminus K$.
    \end{enumerate} 
\end{definition}
\begin{remark}[On the inclusion of these spaces in $\mathcal{L}^p$ spaces]
    When $\mathcal{B}(M)\subset \Sigma_M$, we have the following order of inclusion $\Cc(M,N) \subset \Ccr(M,N)\subset \cLb(M,N)$. Hence, \cref{rem:bounded_measurabe_maps} also shows that, if $h\in \Lb(M,N)$, both sets are also subsets of $\mathcal{L}^\infty_h(M,N)$ for any choice of measure $\mu_M$ and subsets of $\cLph(M,N)$, $p\in [1,\infty)$, when $\mu_M$ is finite. Furthermore, when $h\equiv  z_0\in N$ and $\mu_M$ is infinite but is finite on compact sets, $\Cc(M,N)\cap \cLph(M,N)$ coincides with the set of continuous mappings that are equal to $z_0$ outside a compact set.
\end{remark}

\subsection{On the regularity of measures}

In the following section, we will see that the density of continuous mappings, hence of smooth mappings, heavily relies on the \emph{regularity} of the measure $\mu_M$, so we should clarify what that means. The following definition is based on \cite[342A Definitions (a)]{fremlin2000measure}.

\begin{definition}[Regularity of measures]
Let $(\mathcal{F},\mathcal{A})\in ({2^M})^2$. A measure $\mu_M$ on $(M,\Sigma_M)$ is called:
 \begin{enumerate}[label=(\roman*)]
     \item \emph{inner regular} on $\mathcal{F}$ with respect to $\mathcal{A}$ when each $B\in\mathcal{F}\cap \Sigma_M$ satisfies 
     $\mu_M(B) = \sup \left\{ \mu_M(A): A\in \mathcal{A}\cap \Sigma_M, A\subset B\right\}$, that is, the $\mu_M$-measure of $B$ can be approximated from within by the $\mu_M$-measure of sets in $\mathcal{A}$.
     \item \emph{outer regular} on $\mathcal{F}$ with respect to $\mathcal{A}$ when each $B\in\mathcal{F}\cap \Sigma_M$ satisfies 
     $\mu_M(B) = \inf \left\{ \mu_M(A): A\in \mathcal{A}\cap \Sigma_M, B\subset A\right\}$, that is, the $\mu_M$-measure of $B$ can be approximated from outside by the $\mu_M$-measure of sets in $\mathcal{A}$.

\end{enumerate}
\end{definition}

From this definition, we now clarify the type of regularity that will be used in the proofs of the following sections.

\begin{assumption}[Regularity of $\mu_M$ on a topological base space]
\label{assum:measure_regularity}
From now on, $M$ is further assumed to be a topological space, and we will alternate between the following three regularity assumptions on $\mu_M$:
\begin{enumerate}[
    label=\textsc{(oro)},
    ref=\textsc{(oro)},
]
    \item \label{itm:oro} $\mu_M$ is outer regular on $\mathcal{C}_M\cap \mathcal{F}_{\mu_M}$, the set of closed subsets of $M$ with finite $\mu_M$-measure, with respect to $\mathcal{O}_M$, the set of open subsets of $M$.
\end{enumerate}

\begin{enumerate}[
    label=\textsc{(irc)},
    ref=\textsc{(irc)},
]
    \item \label{itm:irc} $\mu_M$ is inner regular on $\mathcal{F}_{\mu_M}$ with respect to $\mathcal{C}_M$.
\end{enumerate}

\begin{enumerate}[
    label=\textsc{(irk)},
    ref=\textsc{(irk)},
]
    \item \label{itm:irk} $\mu_M$ is inner regular on $\mathcal{F}_{\mu_M}$ with respect to $\mathcal{K}_M$, the set of compact subsets of $M$.
\end{enumerate}
Note that any measure satisfying \ref{itm:irk} also satisfies \ref{itm:irc} as soon as $M$ is Hausdorff.
\end{assumption}
These regularity conditions are standard and easily satisfied. In particular, we are interested in the co-occurrence of \ref{itm:oro} and one of \ref{itm:irc} or \ref{itm:irk}.
First, regarding the co-occurrence of \ref{itm:oro} and \ref{itm:irc}:

\begin{example}[Measures satisfying both \ref{itm:oro} and \ref{itm:irc}] A sufficient condition for $\mu_M$ to satisfy \ref{itm:oro} and \ref{itm:irc} is that it be a finite Borel measure on a topological space $M$ in which open sets are $F_\sigma$ sets, that is, countable unions of closed sets (see \cite[Lemma 7.2.4]{cohn2013measure} for a proof of this result). Metric spaces are prominent examples of such spaces. 
\end{example}

Now, regarding the (usually stronger) co-occurrence of \ref{itm:oro} and \ref{itm:irk}:

\begin{example}[Measures satisfying both \ref{itm:oro} and \ref{itm:irk}]
\label{ex:regular_measures}
    Sufficient conditions to satisfy both \ref{itm:oro} and \ref{itm:irk} are for $\mu_M$ to be regular, in the sense of \cite[p.~189]{cohn2013measure}, on a Hausdorff space $M$ in which each open set is $F_\sigma$ \cite[Proposition~7.2.6]{cohn2013measure}. Examples of regular measures include Radon measures, that is, Borel measures that are finite on compact sets, on a locally compact separable metric space $(M,d_M)$ \cite[Proposition~7.2.3]{cohn2013measure}. This first family of measures includes the Lebesgue measure $\mathscr{L}^d$ on $(\mathbb{R}^d,\lVert \cdot \rVert_d)$ for $d\in \mathbb{N}^*$ and the volume measure on any finite-dimensional Riemannian manifold. The assumption of local compactness on $M$ can be replaced by completeness when considering finite Borel measures on $(M,d_M)$ a complete and separable  (i.e.~Polish) metric space \cite[Proposition~8.1.12]{cohn2013measure}. This second family of measures includes the Wiener measure on $(\mathcal{C}([0,1],\mathbb{R}),d_\infty)$ or probability measures on separable Banach spaces. Also, note that in the first family of measures, \cite[Proposition~7.2.5]{cohn2013measure} imposes the $\sigma$-finiteness of the measure. 
    
    \noindent For a less standard example that is not included in the two previous families of measure spaces, consider $(\mathbb{Q},\lvert \cdot\rvert)$ equipped with the measure $\mu_{\mathbb{Q}}\coloneqq \sum_{n\in\mathbb{N}}2^{-n}\delta_{q_n}$, where $\left\{q_n\right\}_{n\in\mathbb{N}}$ is such that $\mathbb{Q}= \left\{q_n:n\in\mathbb{N}\right\}$ and $\delta_{q_n}$ denotes the Dirac measure at $q_n$ on the Borel $\sigma$-algebra $\mathcal{B}(\mathbb{Q})$, which coincides with the power set of $\mathbb{Q}$. Despite $\mathbb{Q}$ being neither complete nor locally compact, $\mu_{\mathbb{Q}}$ is inner regular on $\mathcal{F}_{\mu_\mathbb{Q}}$ with respect to $\mathcal{K}_{\mathbb{Q}}$. Indeed, note that $\mu_{\mathbb{Q}}$ is finite as $\mu_{\mathbb{Q}}(\mathbb{Q})= \sum_{n\in\mathbb{N}} 2^{-n} < \infty$. Now, let $A\in \mathcal{B}(\mathbb{Q})$. By definition, there exists a subsequence $(n_k)_{k\in\mathbb{N}}$ such that $A=\cup_{k\in\mathbb{N}} \{q_{n_k}\}$. Defining for all $l\in \mathbb{N}$ the sets $K_l \coloneqq \cup_{k\leq l}\{q_{n_k}\}\in \mathcal{B}(\mathbb{Q})$, which are compact as finite unions of compact sets, we get that $(K_l)_{l\in \mathbb{N}}$ is an increasing sequence of sets such that $A=\cup_{l\in\mathbb{N}} K_l$ and $\mu_{\mathbb{Q}}(A) = \lim_{l\to \infty} \mu_{\mathbb{Q}}(K_l)$ \cite[Proposition~1.2.5~(a)]{cohn2013measure}. Hence, $\mu_{\mathbb{Q}}$ is inner regular on $\mathcal{B}(\mathbb{Q})$, thus on $\mathcal{F}_{\mu_{\mathbb{Q}}}$, with respect to $\mathcal{K}_{\mathbb{Q}}$.
\end{example}

We now move on to the approximation of simple mappings by continuous mappings. 

\subsection{Approximation of simple mappings by continuous mappings}

Recall that a topological space $M$ is said to be \emph{normal} if any two disjoint closed sets can be separated by disjoint open sets and \emph{Hausdorff} if any two distinct points can be separated by disjoint open sets.
\begin{theorem}[Approximation of simple mappings by continuous mappings]
\label{th:simple_approx_by_continuous}
Let $p\in [1,\infty)$. Suppose that $M$ is a topological space, that $\Sigma_M$ includes $\mathcal{B}(M)$, that $\mu_M$ satisfies \ref{itm:oro} and that $N$ is path-connected. Then, for any choice of simple mapping $g\in \mathcal{E}(M,N)$ that is equal to $z_0\in N$ outside some $B\in \mathcal{F}_{\mu_M}$, there exists a compact $K_g\subset N$ such that for all $\varepsilon > 0$:
\begin{enumerate}[label=(\roman*)]
    \item \label{cond:irc} if $\mu_M$ satisfies \ref{itm:irc} and $M$ is normal, there exists $g_\varepsilon\in \mathcal{C}(M,N)$, so that $D_p(g,g_\varepsilon) < \varepsilon$, $g_\varepsilon(M)\subset K_g$ and $\overline{g^{-1}_\varepsilon(N\setminus \{z_0\})}\in \mathcal{F}_{\mu_M}$.
    \item \label{cond:irk} if $\mu_M$ satisfies \ref{itm:irk} and $M$ is locally compact and Hausdorff, there exists $g_\varepsilon\in \mathcal{C}(M,N)$, so that $D_p(g,g_\varepsilon) < \varepsilon$, $g_\varepsilon(M)\subset K_g$ and $\overline{g^{-1}_\varepsilon(N\setminus \{z_0\})}\in \mathcal{K}_M\cap \mathcal{F}_{\mu_M}$.
\end{enumerate}

\end{theorem}

\begin{remark}[On the assumption on $g$]
    The condition on $g$ avoids requiring that $\mu_M$ is finite by reducing the approximation problem on $M$ to an approximation on a measurable set of finite $\mu_M$-measure. 
\end{remark}

\begin{remark}[On the assumption of normality]
    Normality is a standard condition in topology satisfied, in particular, by paracompact spaces \cite[Theorem 5.1.5]{engelking1977general}, that is, Hausdorff spaces such that any open cover has a locally finite refinement, which includes metric spaces \cite[Theorem 5.1.3]{engelking1977general}. Also, note that we do not require that every singleton is closed in our definition of normality, nor that the space is additionally Hausdorff, so that the definition we used is sometimes referred to as the \enquote{$T_4$} separation axiom in topology \cite[Section 2, p.~11]{steen1978counterexamples}.
\end{remark}

Before proving \cref{th:simple_approx_by_continuous}, we recall the following well-known separation results for normal spaces.
    \begin{lemma}[Separation of closed sets in normal topological spaces]
\label{lem:separation}
Let $M$ be a normal topological space. Then, the following assertions hold: 
\begin{enumerate}[label=(\roman*)]
\item \label{itm:separation_multiple} for every finite family $(C_i)_{1\leq i\leq n}$ of pairwise disjoint closed sets, there exist pairwise disjoint open sets $(U_i)_{1\leq i\leq n}$ such that $C_i \subset U_i$ for all $i$.  
\item \label{itm:separation} for each closed set $C \subset M$ and open set $U \subset M$ with $C \subset U$, there exists an open set $V$ such that
$C \subset V \subset \overline{V} \subset U.$
\end{enumerate}
\end{lemma}

\begin{proof}
See \ref{appendix:separation}.
\end{proof}

When further assuming that $M$ is locally compact, we get a stronger, also well-known, result for compact sets, and the assumption of normality on $M$ can be replaced with $M$ being Hausdorff. Also, note that a locally compact Hausdorff topological space might not be normal (see the example of the deleted Tychonoff plank \cite[Counterexample 87]{steen1978counterexamples} or of the rational sequence topology \cite[Counterexample 65]{steen1978counterexamples}). 

\begin{lemma}[Separation of compact sets in locally compact Hausdorff topological spaces]
\label{lem:separation_compact}
Let $M$ be a locally compact Hausdorff topological space. Then, the following assertions hold: 
\begin{enumerate}[label=(\roman*)]
\item \label{itm:separation_multiple_compact} for every finite family $(K_i)_{1\leq i\leq n}$ of pairwise disjoint compact sets, there exist pairwise disjoint open sets $(U_i)_{1\leq i\leq n}$ such that $K_i \subset U_i$ for all $i$.  
\item \label{itm:separation_compact} for every compact set $K \subset M$ and open set $U \subset M$ with $K \subset U$, there exists an open set $V$ with compact closure such that
$K \subset V \subset \overline{V} \subset U.$
\end{enumerate}
\end{lemma}
\begin{proof}
    See \ref{appendix:separation_compact}.
\end{proof}

Another key argument in the proof of \cref{th:simple_approx_by_continuous} is the existence of continuous mappings that separate disjoint closed sets. The existence of such mappings in normal topological spaces is given by Urysohn's lemma, whose proof can be found in \cite[Theorem 1.5.11]{engelking1977general} or \cite[Chapter II,  \S 4, Theorem 4.2]{lang2012real}. Note that in both references, \enquote{normality} means normality in the sense defined above \emph{and} Hausdorff, but the Hausdorff assumption is unnecessary in their proofs.

\begin{theorem}[Urysohn's lemma in normal topological spaces]
\label{th:urysohn}
    Let $M$ be a normal topological space and let $C$ and $F$ be two disjoint closed sets. Then, there exists a continuous mapping $I:M\to [0,1]$ such that $I|_C\equiv 1$ and $I|_F\equiv 0$.
\end{theorem}

It is then well-known that, similarly to the two previous lemmas, normality can be replaced in Urysohn's lemma by instead assuming that $M$ is locally compact and Hausdorff to obtain a continuous mapping separating two disjoint sets, one being compact and the other being closed. 

\begin{corollary}[Urysohn's lemma in locally compact Hausdorff topological spaces]
\label{cor:urysohn}
    Let $M$ be a locally compact Hausdorff topological space and let $K$ and $F$ be two disjoint sets with $K$ compact and $F$ closed. Then, there exists a continuous mapping $I:M\to [0,1]$ such that $I|_K\equiv 1$ and $I|_F\equiv 0$.
\end{corollary}

\begin{proof}
    See \ref{appendix:ursyohn}.
\end{proof}

\begin{proof}[Proof of \cref{th:simple_approx_by_continuous}]
    Let $p\in [1,\infty)$, $\varepsilon > 0$ and $g \in \mathcal{E}(M,N)$ such that $g_{M\setminus B}\equiv z_0\in N$ for some $B\in \mathcal{F}_{\mu_M}$. Also, we can assume that $z_0\notin g(B)$. 

    \begin{enumerate}[label=(\roman*),wide]
        \item We can assume that there exists a finite indexing set $I$ such that $g(B) = \left\{y_i\in N: i\in I\right\}$ and define $B_i \coloneqq g^{-1}(\left\{y_i\right\})$. The proof then essentially consists in finding a \enquote{sufficiently large} closed $C_i\subset B_i$ and a \enquote{sufficiently small} open $V_i \supset C_i$ such that the approximation can be constructed via a continuous interpolation on $V_i\setminus C_i$, while managing potential overlaps to deal with the nonlinear nature of $N$.

        \par\medskip \noindent \emph{Step 1 (Approximation by a modification of $g$ that is constant outside \enquote{sufficiently large} closed sets):} Since $\mu_M$ satisfies \ref{itm:irc}, there exists a family of closed sets $(C_i)_{i\in I}$ such that $C_i\subset B_i$ and $\mu_M(B_i \setminus C_i) < (\varepsilon/(2\lvert I\rvert^{1/p} R_0))^p$, with $R_0\coloneqq \max_{i\in I} d_N(y_i,z_0)$. Thus, defining the mapping $\tilde{g}: M\to N$ as
    $$\tilde{g}(x)\coloneqq \begin{cases}
        g(x)&\text{if $x\in \cup_i C_i$}\\
        z_0&\text{otherwise}
    \end{cases}$$
    yields that, by construction, $D_p(g,\tilde{g})^p \leq R_0^p\,\sum_{i\in I} \mu_M(B_i \setminus C_i) < (\varepsilon/2)^p$.
    
    \par\medskip \noindent \emph{Step 2 (Construction of a continuous approximation by continuously interpolating between the \enquote{sufficiently large} closed sets and \enquote{sufficiently small} open sets):} Since $N$ is path-connected, there exists for all $i\in I$ a continuous curve $\gamma_i: [0,1]\to N$ such that $\gamma_i(0) = z_0$ and $\gamma_i(1)=y_i$. In addition, define $K_g \coloneqq \cup_{i\in I} \gamma_i([0,1])$, which is compact as the finite union of compact sets. From \ref{itm:oro}, there exists a family of open sets $(U_i)_{i\in I}$ such that $C_i\subset U_i$ and $\mu_M(U_i\setminus C_i) < (\varepsilon/(2\lvert I\rvert^{1/p} R_1))^p$, with $R_1 \coloneqq \max_{(i,t)\in I\times [0,1]} d_N(\gamma_i(t),z_0)$. Using \ref{itm:separation_multiple} of \cref{lem:separation}, we can assume that the $U_i$'s are disjoint. Then, using \ref{itm:separation} of \cref{lem:separation}, we get the existence of a family $(V_i)_{i\in I}$ such that $C_i\subset V_i\subset \overline{V_i}\subset U_i$. Then, recalling that for all $i\in I$ the closed sets $C_i$ and $M\setminus V_i$ are disjoint, we have, by Urysohn's lemma (\cref{th:urysohn}), that there exists a continuous mapping $I_i: M \to [0,1]$ such that $I_i|_{C_i} \equiv 1$ and $I_i|_{M\setminus V_i} \equiv 0$. We thus have, by composition, that the mapping $\gamma_i \circ I_i$ is continuous on $M$ and constant outside $V_i$. Therefore, define the mapping $g_\varepsilon : M\to N$ as $g_\varepsilon|_{U_i} \coloneqq \gamma_i\circ I_i|_{U_i}$ for all $i\in I$ and $g_\varepsilon|_{M\setminus U} \equiv z_0$ with $U\coloneqq \cup_{i\in I} U_i$. By construction, $g_\varepsilon$ is continuous. Indeed, let $x\in M$. If $x\in U$, there exists $i\in I$ such that $x\in U_i$ and we know that $g_\varepsilon|_{U_i} = \gamma_i \circ I_i|_{U_i}$. Otherwise, $x\in M\setminus U$, so that, if $V\coloneqq \cup_{i\in I} V_i$, $x\in M\setminus \overline{V}$. Hence, there exists a neighborhood $U_x\subset M\setminus\overline{V}$ such that $g_\varepsilon|_{U_x} \equiv z_0$. In both cases, there exists a neighborhood $U_x$ of $x$ such that $g_\varepsilon|_{U_x}$ is continuous.  
    Also, note that $g_\varepsilon(M)$ is a subset of the compact $K_g$. By construction of $g_\varepsilon$, we have that
    $$D_p(\tilde{g},g_\varepsilon)^p = \sum_{i\in I}\int_{V_i\setminus C_i} d_N(z_0,\gamma_i\circ I_i(x))^p\dif\mu_M(x)\leq R_1^p \sum_{i\in I} \mu_M(V_i\setminus C_i)\leq R_1^p \sum_{i\in I} \mu_M(U_i\setminus C_i)< (\varepsilon/2)^p.$$

    \noindent Finally, we get, by the triangle inequality, that $D_p(g,g_\varepsilon) \leq D_p(g,\tilde{g})+ D_p(\tilde{g},g_\varepsilon) < \varepsilon$.
    \item The proof works in the exact same way by taking compact $C_i$'s using the \ref{itm:irk} regularity of $\mu_M$, taking $U_i$'s still using \ref{itm:oro} but using \ref{itm:separation_multiple_compact} of \cref{lem:separation_compact} to make them disjoint, taking $V_i$'s using \ref{itm:separation_compact} of \cref{lem:separation_compact} and using a corollary of Urysohn's lemma holding in locally compact Hausdorff topological spaces to pick the continuous transition mappings (\cref{cor:urysohn}). Hence, $g_\varepsilon$ is constant outside the compact $\cup_{i\in I} \overline{V}_i$, as the finite union of compact sets.\qedhere
    \end{enumerate}
\end{proof}

As in \cref{sec:density_simple}, we now discuss sharpness of the assumptions used in \cref{th:simple_approx_by_continuous} by providing counterexamples. The reader may skip to \cref{sec:extension_lebesgue} in the first reading.
Beginning with the outer regularity condition on $\mu_M$ required in \ref{cond:irc} of \cref{th:simple_approx_by_continuous}, we provide a counterexample when this condition is not satisfied.

\begin{cexample}[When $\mu_M$ does not satisfy \ref{itm:oro}]
    Suppose that $M=\mathbb{R}$ with its standard topology, that $N=\mathbb{R}_+$ with the standard metric on $\mathbb{R}$ and define the Borel measure $\mu_M = \rho\mathscr{L}^1 + \delta_0$ with $\rho(x)\coloneqq \frac{1}{\lvert x\rvert}$ for all $x\in \mathbb{R}^*$ and $\rho(0) \coloneqq 0$. Then, $\mu_M$ does not satisfy \ref{itm:oro} since, for every open subset $U$ of $M$ containing $0$, $\mu_M(U) = \infty$, but $\mu_M(\left\{0\right\}) = 1$. On the other hand, $\mu_M$ satisfies \ref{itm:irk}, hence \ref{itm:irc}, since $\rho\mathscr{L}^1$ is a Radon measure on $\mathbb{R}\setminus \{0\}$ and $\delta_0$ is concentrated on the compact set $\{0\}$. Thus, the mapping $g\coloneqq \mathds{1}_{\{0\}}$ belongs to $\mathcal{E}(\mathbb{R},\mathbb{R}_+)$, but $g$ cannot be approximated arbitrarily closely by continuous mappings for $D_p$, regardless of the value of $p\in [1,\infty)$. Indeed, let $f\in \mathcal{C}(\mathbb{R},\mathbb{R}_+)$. If $f(0) > 0$, there exists a neighborhood $U_0$ of $0$ such that $f(x)\geq f(0)/2$ for every $x\in U_0$. Therefore, $D_p(g,f)^p \geq (f(0)/2)^p \mu_M(U_0\setminus \{0\}) = \infty$. If $f(0) = 0$, then $D_p(g,f)^p \geq \lvert g(0) - f(0)\rvert^p\mu_M(\{0\}) = 1$. In both cases, we get $D_p(g,f) \geq 1$.
\end{cexample}

Now, we provide two counterexamples when $g$ is simple but not constant outside a set of finite $\mu_M$-measure.
\begin{cexample}[When $g$ is not constant outside a measurable set of finite $\mu_M$-measure]
    Let us provide the following two counterexamples:
    \begin{enumerate}[label=(\roman*),wide]
        \item Suppose that $N=\mathbb{R}$ with its standard metric and that $M=\{0\}\cup \{1/n: n\in\mathbb{N}^*\}$ with the topology induced by the standard topology of $\mathbb{R}$, so that $M$ is Hausdorff and compact, hence normal \cite[Theorem 3.1.9]{engelking1977general}. Let $\mu_M$ be the Borel measure such that $\mu_M(A)= \infty$ if $0\in A$ and $\mu_M(A)= \#(A)$ otherwise, with $\#$ denoting the counting measure. Then, $\mu_M$ satisfies both \ref{itm:oro} and \ref{itm:irk}, hence also \ref{itm:irc}. Indeed, every set of finite $\mu_M$-measure is a finite subset of $M\setminus\{0\}$, hence is compact and closed in $M$, while every subset of $M\setminus \{0\}$ is open in $M$. Now, consider the mapping $g\coloneqq \mathds{1}_{\{0\}}$ that is simple and is not constant outside a measurable set of finite $\mu_M$-measure. However, there is no continuous mapping approximating $g$ for $D_p$, regardless of the value of $p\in [1,\infty)$. Indeed, let $f\in \mathcal{C}(M,\mathbb{R})$ and suppose that $D_p(f,g) <\infty$. Since $\mu_M(\{0\})= \infty$, we must have $f(0)=g(0)=1$. Thus, by continuity of $f$, there exists a neighborhood $U_0$ of $0$ in $M$ such that $\lvert f(x)\rvert\geq 2^{-1}$ for every $x\in U_0$. Since, $g(x)=0$ for every $x\in U_0\setminus\{0\}$ and $\mu_M(U_0\setminus \{0\}) =\infty$ (since $U_0\setminus\{0\}$ is infinite), we get that $$D_p(f,g)^p\geq \int_{U_0\setminus \{0\}} \lvert f(x)\rvert^p\dif \mu_M(x)\geq 2^{-p} \mu_M(U_0\setminus \{0\}) = \infty,$$
        which yields a contradiction.
        \item Suppose that $M= N =\mathbb{R}$ with its standard topology, and that $\mu_M$ is the purely infinite Borel measure such that $\mu_M(A)=\infty$ on any nonempty Borel set. Then, $\mu_M$ satisfies both \ref{itm:oro} (since $\mathcal{F}_{\mu_M}=\{\emptyset\}$) and \ref{itm:irk} (since the empty set is compact). Now, consider the mapping $g\coloneqq \mathds{1}_{\mathbb{R}_+}$ that is simple and is not constant outside a measurable set of finite $\mu_M$-measure.  However, there is no continuous mapping approximating $g$ for $D_p$, regardless of the value of $p\in [1,\infty)$. Indeed, let $f\in \mathcal{C}(\mathbb{R},\mathbb{R})$ and suppose that $D_p(f,g) <\infty$. Then, we have $f=g$, which yields a contradiction since $g$ is discontinuous.
    \end{enumerate}
\end{cexample}

Also, we provide a counterexample when the target space is not path-connected.

\begin{cexample}[When $N$ is not path-connected]
Suppose $M= \mathbb{R}$ is equipped with the $\sigma$-algebra $\Sigma_M= \mathcal{B}(\mathbb{R})$ and the measure $\mu_M=\mathscr{L}^1$, which satisfies both \ref{itm:oro} and \ref{itm:irk}, and that $N =\mathbb{R}^*$ is equipped with its standard metric. Then consider $g\coloneqq -\mathds{1}_{[-1,1]} + \mathds{1}_{\mathbb{R}\setminus [-1,1]}$ which belongs to $\mathcal{E}(\mathbb{R},\mathbb{R}^*)$ and is constant outside the Borel measurable set $[-1,1]$ of finite Lebesgue measure. Then, $g$ cannot be approximated by a continuous mapping $f$ which would necessarily have its values included only in one of the two connected components of $\mathbb{R}^*$.
\end{cexample}

Now, we provide a counterexample to \ref{cond:irk} of \cref{th:simple_approx_by_continuous} when $\mu_M$ satisfies \ref{itm:irc} but not \ref{itm:irk}.

\begin{cexample}[When $\mu_M$ does not satisfy \ref{itm:irk} but satisfies \ref{itm:irc}]
Let $p\in [1,\infty)$. Denote by $\omega_1$ the first uncountable ordinal and set $M=[0,\omega_1)$, the set of countable ordinals, equipped with the order topology \cite[4A2R,~4A2S]{fremlin2003topological}, so that $[0,\omega_1)$ is Hausdorff and locally compact, and with $\Sigma_M= \mathcal{B}([0,\omega_1))$. In addition, let $\mu_M$ be Dieudonné's measure \cite[411Q,~4A3J]{fremlin2003topological}, that is, define $\mu_M$ as  
$$\mu_M(A)\coloneqq \begin{cases}
    1,&\text{if $A$ contains a closed cofinal subset of $[0,\omega_1)$,}\\
    0,&\text{if $A$ is disjoint from a closed cofinal subset of $[0,\omega_1)$},
\end{cases},$$
where we recall that a set $A\subset [0,\omega_1)$ is called cofinal if for every $x\in [0,\omega_1)$ there exists $x'\in A$ such that $x\leq x'$, and is closed if $\sup B\in A$ whenever $\emptyset\neq B\subset A$ and $\sup B< \omega_1$ \cite[4A1B (a)]{fremlin2003topological}. Also, suppose that $N=\mathbb{R}$, equipped with its standard metric. Then, $\mu_M$ satisfies both \ref{itm:oro} and \ref{itm:irc} \cite[411Q]{fremlin2003topological}, but does not satisfy \ref{itm:irk}. Indeed, if $K$ is compact in $[0,\omega_1)$, it is bounded \cite[4A2R (g)]{fremlin2003topological}, so $K\subset [0,\beta]$ for some $\beta < \omega_1$ . Therefore, $K$ and $[\beta +1, \omega_1)$ are disjoint and $[\beta +1, \omega_1)$ is a closed and cofinal subset of $[0,\omega_1)$ \cite[4A1B]{fremlin2003topological}, hence $\mu_M(K)=0$. In particular, if $C$ is a closed cofinal subset of $[0,\omega_1)$, we have $\mu_M(C)=1 > \sup\{\mu_M(K): K\in\mathcal{K}_{[0,\omega_1)},K\subset C\}=0$. Now, define $g\coloneqq \mathds{1}_{C}$.
Then, let $f\in \mathcal{C}([0,\omega_1),\mathbb{R})$ such that $\operatorname{supp}(f)\coloneqq \overline{\{x\in [0,\omega_1): f(x)\neq 0\}}$ is compact in $[0,\omega_1)$. Since $\operatorname{supp}(f)$ is compact, it has null $\mu_M$-measure. Hence, we have 
$$D_p(g,f)^p\geq \int_{C\setminus \operatorname{supp}(f)}\lvert g(x)- f(x)\rvert^p\dif \mu_M(x) =\mu_M(C\setminus \operatorname{supp}(f))= \mu_M(C) =1 > 0.$$
Hence, there is no $f\in \mathcal{C}([0,\omega_1),\mathbb{R})$ such that $\operatorname{supp}(f)$ is compact and $D_p(f,g) < 1$. Thus, the conclusion of assertion  \ref{cond:irk} in \cref{th:simple_approx_by_continuous} fails for $z_0= 0$.
\end{cexample}

Finally, we provide a counterexample to \ref{cond:irk} of \cref{th:simple_approx_by_continuous} when $\mu_M$ satisfies both \ref{itm:oro} and \ref{itm:irk}, but the base space is not locally compact.

\begin{cexample}[When $M$ is not locally compact]
  Let $M=\mathbb{H}$ be an infinite-dimensional separable Hilbert space equipped with $\Sigma_M=\mathcal{B}(\mathbb{H})$ and a nontrivial finite Borel measure $\mu_M$, and let $N=\mathbb{R}$ equipped with its standard metric. Since $\mathbb{H}$ is a Polish space and $\mu_M$ is a finite Borel measure, $\mu_M$ satisfies both \ref{itm:oro} and \ref{itm:irk} by \cref{ex:regular_measures}. The only continuous mapping $f:\mathbb{H}\to\mathbb{R}$ with compact support is the null mapping. Indeed, if $f(x_0)\neq 0$ for some $x_0\in \mathbb{H}$, then continuity gives the existence of $r>0$ such that the closed ball $\overline{B}(x_0,r)$ is included in the (compact) support of $f$. Hence, as a continuous image of $\overline{B}(x_0,r)$, the unit closed ball would be compact, contradicting Riesz's theorem \cite[Chapter II,  \S 3, Corollary 3.15]{lang2012real}. Therefore, a simple mapping that is non-null on a set of positive $\mu_M$-measure cannot be approximated in $D_p$ by continuous mappings with compact support, for any $p\in[1,\infty)$.
\end{cexample}

We now move on to the extension of the previous approximation result to Lebesgue mappings. 

\subsection{Extension to Lebesgue mappings}
\label{sec:extension_lebesgue}

Together with \cref{th:density_simple}, we retrieve, just as in the linear case, that Lebesgue mappings can be approximated by continuous mappings.

\begin{theorem}[Density of continuous mappings]
\label{th:density_continuous}
Let $h\in \Delta(M,N)$ (\cref{def:constant_mappings}) and $p\in [1,\infty)$. Suppose that $M$ is a topological space, that $\Sigma_M$ includes $\mathcal{B}(M)$, that $\mu_M$ satisfies \ref{itm:oro} and that $N$ is path-connected. Then,
\begin{enumerate}[label=(\roman*),wide]
    \item if $\mu_M$ satisfies \ref{itm:irc} and $M$ is normal, $\Ccr(M,N)\cap \cLph(M,N)$ forms a dense subspace of $\cLph(M,N)$ and so determines a dense subspace of $\Lph(M,N)$.
    \item if $\mu_M$ satisfies \ref{itm:irk} and $M$ is locally compact and Hausdorff, $\Cc(M,N)\cap \cLph(M,N)$ forms a dense subspace of $\cLph(M,N)$ and so determines a dense subspace of $\Lph(M,N)$.
\end{enumerate}
\end{theorem}

\begin{remark}[Relaxation of the assumption on $h$ when $\mu_M$ is finite]
    When $\mu_M$ is finite, \cref{th:density_continuous} holds for any choice of bounded base mapping $h\in \Lb(M,N)$, by \cref{prop:indep_def}.
\end{remark}

\begin{proof}[Proof of \cref{th:density_continuous}]
Let $h\equiv z_0\in N$, $p\in[1,\infty)$, $f\in \cLph(M,N)$ and $\varepsilon > 0$. 

\begin{enumerate}[label=(\roman*),wide]
    \item 
    \emph{Step 1 (Approximation by a simple mapping):} By \cref{th:density_simple}, there exists $g\in \mathcal{E}(M,N)\cap \cLph(M,N)$ such that $D_p(f,g)\leq \varepsilon/2$. Since $h\equiv z_0$, $\mu_M(g^{-1}(N\setminus \{z_0\}))<\infty$.
    
\par\medskip \noindent \emph{Step 2 (Approximation by a continuous mapping):} By \ref{cond:irc} of \cref{th:simple_approx_by_continuous}, there exists $g_\varepsilon \in \Ccr(M,N)$ such that $D_p(g,g_\varepsilon) < \varepsilon/2 $, hence $g_\varepsilon\in \Ccr(M,N)\cap \cLph(M,N)$. 

\noindent Finally, we get, by the triangle inequality, that $D_p(f,g_\varepsilon) \leq D_p(f,g)+D_p(g,g_\varepsilon) < \varepsilon$.
\item The proof follows the same line, but picking $g_\varepsilon\in \Cc(M,N)$ thanks to \ref{cond:irk} of \cref{th:simple_approx_by_continuous}.\qedhere
\end{enumerate}
\end{proof}

\begin{remark}[Related results in the literature] A similar result was proved by S.~Neumayer, J.~Persch and G.~Steidl \cite[Theorem~2.2]{neumayer2018} in the case where $h$ is constant, $M$ is an open, bounded, and connected (Lipschitz) domain of $\mathbb{R}^n$, with $n\geq 2$, $\mu_M = \mathscr{L}^n$ and $N$ is a locally compact Hadamard space. Since the Lebesgue measure satisfies \ref{itm:oro} and \ref{itm:irk} (\cref{ex:regular_measures}) and by the same arguments exposed in \cref{rem:cor_neumayer}, their result comes as a corollary of \cref{th:density_continuous}.
\end{remark}

\section{Density of smooth mappings}
\label{sec:density_smooth}

\cref{th:density_continuous} extends to smooth mappings (\cref{def:smooth_mappings}) when the base space is a smooth Banach manifold (\cref{def:banach_manifolds}) such that every open cover has a smooth partition of unity subordinated to it (\cref{def:partition_unity}) and the target space is a smooth and connected Banach manifold equipped with a metric whose induced topology is weaker than or equal to the manifold topology. This extension essentially follows from the same arguments as in the proofs of \cref{th:simple_approx_by_continuous} and \cref{th:density_continuous} by additionally using smooth paths and approximating the continuous transition mappings by smooth mappings using smooth partitions of unity, as done in the proof of Whitney's approximation theorem \cite{whitney1936differentiable} provided in \cite[Theorem 6.21]{lee2003smooth}. To some extent, the regularization using smooth partitions of unity acts as a surrogate to regularization by convolution with mollifiers used to prove the density of smooth mappings in the linear case. In what follows, given $r\in \mathbb{N}^*\cup\{\infty\}$, we use the notation 
$$\Crcr(M,N)\coloneqq \mathcal{C}^r(M,N)\cap \Ccr(M,N),\quad \Crc(M,N)\coloneqq \mathcal{C}^r(M,N)\cap \Cc(M,N).$$

 \begin{theorem}[Density of smooth mappings]
 \label{th:density_smooth}
        Let $h\in \Delta(M,N)$ (\cref{def:constant_mappings}), $p\in [1,\infty)$ and $r\in \mathbb{N}^*\cup \{\infty\}$.
        Suppose that $M$ is a $\mathcal{C}^r$ Banach manifold (not necessarily Hausdorff) such that every open cover has a $\mathcal{C}^r$ partition of unity subordinated to it, that $\Sigma_M$ includes $\mathcal{B}(M)$, that $\mu_M$ satisfies \ref{itm:oro}, that $N$ is a connected $\mathcal{C}^r$ Banach manifold and that $d_N$ is a metric on $N$ whose induced topology is weaker than (or equal to) the manifold topology. Then,
        \begin{enumerate}[label=(\roman*)]
            \item \label{cond:irc_smooth} if $\mu_M$ satisfies \ref{itm:irc} and $M$ is normal, $\Crcr(M,N)\cap \cLph(M,N)$ is a dense subspace of $\cLph(M,N)$ and so determines a dense subspace of $\Lph(M,N)$.
            \item \label{cond:irk_smooth} if $\mu_M$ satisfies \ref{itm:irk} and $M$ is locally compact and Hausdorff, $\Crc(M,N)\cap \cLph(M,N)$ is a dense subspace of $\cLph(M,N)$ and so determines a dense subspace of $\Lph(M,N)$.
        \end{enumerate}
    \end{theorem}

    \begin{remark}[On the assumption of existence of $\mathcal{C}^r$ partitions of unity]
        The existence of $\mathcal{C}^r$ partitions of unity subordinated to every open cover of $M$, for $r\in \mathbb{N}^*\cup \{\infty\}$, follows, for instance, under the following assumptions: $M$ is Lindelöf and Hausdorff, and is modeled on a Banach space $\mathbb{B}$ whose norm is $C^r$ on $\mathbb{B}\setminus \{0_{\mathbb{B}}\}$ (see \cite[Theorem 1]{bonic1966smooth} for the proof in the case where $M$ is separable and \cite[Section 3.3]{lloyd1974smooth} for the extension to the case where $M$ is Lindelöf). Recall that $M$ is Lindelöf if every open cover of $M$ has a countable subcover.
    \end{remark}

    \begin{remark}[Extension beyond Banach manifolds]
        \cref{th:density_smooth} could certainly be extended to manifolds modeled on convenient vector spaces (see \cite[Section~27]{kriegl1997convenient} for a definition) under appropriate conditions ensuring the existence of subordinated smooth partitions of unity for every open cover of $M$ such as \cite[Theorem 16.10]{kriegl1997convenient}, but we do not further discuss this matter here.
    \end{remark}
\begin{proof}[Proof of \cref{th:density_smooth}]
        
\noindent Let $p\in [1,\infty)$, $f\in \cLph(M,N)$ and $\varepsilon > 0$.

\begin{enumerate}[label=(\roman*),wide]
    \item \emph{Step 1 (Approximation by a simple mapping):} By \cref{th:density_simple}, there exists $g\in \mathcal{E}(M,N)\cap \cLph(M,N)$ such that $D_p(f,g)< \varepsilon/3$. Therefore, defining $B\coloneqq g^{-1}(N\setminus \{z_0\})$ and recalling $h\equiv z_0$, 
    we have $\mu_M(B)<\infty$. In addition, we can assume that there exists a finite indexing set $I$ such that $g(B) = \left\{y_i\in N: i\in I\right\}$ and denote $B_i \coloneqq g^{-1}(\left\{y_i\right\})\subset B$.

\par\medskip
  \noindent \emph{Step 2 (Approximation by a continuous map using smooth paths):} Since $N$ is a connected $\mathcal{C}^r$ Banach manifold, for every $i\in I$ there exists a $\mathcal{C}^r$ path $\gamma_i: [0,1]\to N$ such that $\gamma_i(0) = z_0$ and $\gamma_i(1) = y_i$ \cite[\S 3, Lemma 3.1]{palais1966lusternik}. Then, following the proof of \cref{th:simple_approx_by_continuous}, there exist two families of open sets $(V_i)_{i\in I}$ and $(U_i)_{i\in I}$ such that $\overline{V_i} \subset U_i$ with $V_i \in \mathcal{F}_{\mu_M}$, $U_i \cap U_{j} = \emptyset$ whenever $j\neq i$. We also obtain a family $(I_i)_{i\in I}$ of continuous transition mappings  $I_i: M \to [0,1]$ vanishing outside $V_i$ together with a map $g_\varepsilon \in \Ccr(M,N)$ such that $D_p(g,g_\varepsilon) < \varepsilon/3 $. Moreover, $g_\varepsilon|_{U_i} = \gamma_i \circ I_i$ for all $i\in I$ and $g_\varepsilon|_{M\setminus U}\equiv z_0$ with $U\coloneqq \cup_{i\in I}U_i$. In addition, for every $i\in I$ the path $\gamma_i$ is continuous from $[0,1]$ to $N$ endowed with the manifold topology. Since $[0,1]$ is compact, $\gamma_i$ is uniformly continuous from $[0,1]$ to $(N,d_N)$. Consequently, there exists $\delta_i > 0$ such that for every $(s,t)\in [0,1]^2$ satisfying $\lvert s - t\rvert < \delta_i$, we have $d_N(\gamma_i(s),\gamma_i(t)) < \varepsilon/(3^p\lvert I\rvert (1+\mu_M(V_i)))^{1/p}$.

\par\medskip
\noindent \emph{Step 3 (Smoothing of the continuous transition mappings):} 
Fix $i\in I$ and define the set $U_0\coloneqq M\setminus \overline{V_i}$. Then, for every $x\in U_i$ choose, by continuity of $I_i$, an open neighborhood $U_x\subset U_i$ such that $\lvert I_i(x) - I_i(x')\rvert < \delta_i$ for every $x'\in U_x$. We can also assume that $U_x\subset V_i$ whenever $x\in V_i$ since $V_i$ is open. Let $S\coloneqq\{0\}\cup U_i$. Since $\overline{V_i}\subset U_i$, the family $\{U_s\}_{s\in S}$ is an open cover of $M$ and, since every open cover of $M$ has a ($\mathcal{C}^r$) partition of unity subordinated to it (\cref{def:partition_unity}), $\{U_s\}_{s\in S}$ admits a locally finite open refinement $\{W_s\}_{s\in S}$ such that $W_s\subset U_s$ for every $s\in S$ \cite[Lemma~5.1.8]{engelking1977general}. Now, let $\{\rho_s\}_{s\in S}$ be a $\mathcal{C}^r$ partition of unity satisfying $\operatorname{carr}(\rho_s)\coloneqq \{x\in M: \rho_s(x) \neq 0\}\subset W_s$ for every $s\in S$. Since $\{W_s\}_{s\in S}$ is locally finite, so is $\{\rho_s\}_{s\in S}$. Define the map $\tilde{I}_i: M\to \mathbb{R}$ by $\tilde{I}_i(x)\coloneqq \sum_{s\in U_i} \rho_s(x) I_i(s)$. By construction, $\tilde{I}_i$ belongs to $\mathcal{C}^r(M,[0,1])$ as the locally finite convex combination of mappings in $\mathcal{C}^r(M,[0,1])$. Moreover, for every $x\in M$, $\tilde{I}_i$ satisfies that 
\begin{align*}
  \lvert \tilde{I}_i(x) - I_i(x)\rvert &= \bigg\lvert \sum_{s\in U_i} \rho_s(x)I_i(s) - \bigg(\rho_0(x)+ \sum_{s\in U_i} \rho_s(x)\bigg) I_i(x)\bigg\rvert
  \leq \sum_{s\in U_i}\rho_s(x)\lvert I_i(s)- I_i(x)\rvert< \delta_i,
\end{align*}
where we have used the fact that $\rho_0I_i\equiv 0$ since $I_i$ vanishes on $\operatorname{carr}(\rho_0)$. Note also that for $x\in M\setminus V_i$, $\rho_s(x)I_i(s)=0$ since $I_i(s)=0$ for $s\notin V_i$ and $\rho_s(x)=0$ for $s\in V_i$ since $\operatorname{carr}(\rho_s)\subset W_s\subset U_s\subset V_i$. Consequently $\gamma_i \circ \tilde{I}_i$ is identically equal to $z_0$ outside $V_i$. Moreover, by composition, it belongs to $\mathcal{C}^r(M,N)$.

\par\medskip
\noindent \emph{Step 4 (Approximation by a smooth map using the smoothed transition mappings):}
Define $\tilde{g}_\varepsilon: M\to N$ by $\tilde{g}_\varepsilon|_{U_i} = \gamma_i\circ\tilde{I}_i|_{U_i}$ for every $i\in I$ and $\tilde{g}_\varepsilon|_{M\setminus U} \equiv z_0$. This definition is unambiguous because the sets $U_i$, $i\in I$, are pairwise disjoint. Let $V\coloneqq\cup_{j\in I}V_j$. The mapping $\tilde{g}_\varepsilon$ satisfies $\tilde{g}_\varepsilon\equiv z_0$ on $M\setminus V$. Indeed, by Step 3, for every $j\in I$, the map $\gamma_j\circ\tilde{I}_j$  is identically equal to $z_0$ on $M\setminus V\subset  M\setminus V_j$. In particular, $\tilde{g}_\varepsilon$ coincides with $\gamma_i\circ\tilde{I}_i$ on the open set $U_i\cup (M\setminus \overline{V})$.
Since $\gamma_i\circ\tilde{I}_i\in \mathcal{C}^r(M, N)$ for any $i$, it follows that $\tilde{g}_\varepsilon$ is of class $\mathcal C^r$ on $M=\cup_{i\in I} (U_i\cup  (M\setminus \overline{V}))$. In addition, $\tilde{g}_\varepsilon(M)$ is contained in the compact set $ \cup_{i\in I} \gamma_i([0,1])$, so that $\tilde{g}_\varepsilon \in \Crcr(M,N)$. Furthermore, we have
        \begin{align*}
        D_p(g_\varepsilon,\tilde{g}_\varepsilon)^p = \sum_{i\in I}\int_{V_i}d_N(\gamma_i(I_i(x)),\gamma_i(\tilde{I}_i(x)))^p\dif \mu_M(x)< (\varepsilon/3)^p\,.
        \end{align*}

        \noindent Finally, the triangle inequality yields  $D_p(f,\tilde{g}_\varepsilon) \leq D_p(f,g)+ D_p(g,g_\varepsilon) + D_p(g_\varepsilon, \tilde{g}_\varepsilon) < \varepsilon$.

\item The proof follows the same line, but using \ref{cond:irk} of \cref{th:simple_approx_by_continuous} to take $g_\varepsilon\in \Cc(M,N)$.\qedhere 
\end{enumerate}
\end{proof}

\section{Notation and terminology}
\label{sec:notations}

This section clarifies some recurrent notions and notation used in the present article:
\subsection*{Spaces}
\begin{itemize}
    \item $\mathbb{N}$ denotes the set of nonnegative integers (including $0$). 
    \item $\mathbb{N}^*$ denotes the set of positive integers (excluding $0$). 
    \item $\mathbb{R}$ denotes the set of real numbers, also referred to as the \emph{real line}. 
    \item $\mathbb{R}_+$ denotes the set of nonnegative real numbers (including $0$).
    \item $\mathbb{R}^*$ denotes the set of nonzero real numbers (excluding $0$). 
     \item $\mathbb{Q}$ denotes the set of rational numbers.
\end{itemize}
\subsection*{Topology}

\begin{itemize} 
   \item $\mathcal{O}_X$ denotes the set of open subsets of the topological space $X$.
   \item $\mathcal{C}_X$ denotes the set of closed subsets of the topological space $X$.
   \item $\mathcal{K}_X$ denotes the set of compact subsets of the topological space $X$.
    \item $f(X)\coloneqq \left\{ f(x): x\in X\right\}$ denotes the \emph{range} of a mapping $f: X\to Y$ between sets $X$ and $Y$.
    \item $\operatorname{carr}(f)\coloneqq\{x\in X: f(x)\neq 0\}$ denotes the \emph{carrier} of a mapping $f: X\to \mathbb{R}$. 
    \item $\overline{A}$ denotes the \emph{closure} of the subset $A \subset X$ in the topological space $X$, that is, the intersection of all closed subsets of $X$ that include $A$.
    \item A set $X$ is called \emph{countable} if there exists an injective mapping $\varphi: X\to \mathbb{N}$.
    \item A subset $S$ of a (semi-)metric space $(X,d_X)$ is called \emph{dense} when for each element $x$ of $X$ and real number $\varepsilon > 0$ there exists an element $s$ of $S$ such that $d_X(x,s) < \varepsilon$.
    \item A metric space $(X,d_X)$ is called \emph{separable} when it has a countable dense subset.
    \item A sequence $(x_n)_{n\in\mathbb{N}}$ in a metric space $(X,d_X)$ is called \emph{Cauchy} when for all real $\varepsilon > 0$ there exists an integer $n_0$ such that for all $n\geq n_0$ and $m\geq n_0$, $d_X(x_n,x_m) < \varepsilon$.
    \item A metric space $(X,d_X)$ is called \emph{complete} when every Cauchy sequence converges to an element of $X$.
    \item $A^c$ denotes the \emph{complement} of a subset $A$ of a set $X$, that is, $A^c \coloneqq \left\{x \in X: x\notin A\right\}$.
    \item $A\setminus B\coloneqq A \cap B^c$ denotes the \emph{set difference} between two subsets $A$ and $B$ of a set $X$. 
    \item $A\Delta A'\coloneqq (A\setminus A')\cup (A'\setminus A)$ denotes the \emph{symmetric set difference} between two subsets $A$ and $A'$ of a set $X$.
\end{itemize}

\subsection*{Measure theory}
\begin{itemize}
    \item $\mathscr{L}^d$ denotes the $d$-dimensional Lebesgue measure.
    \item $\delta_x$ denotes the Dirac measure at a point $x$ of a set $M$.
    \item $\#$ denotes the counting measure.
    \item A measure $\mu_M$ on a measurable space $(M,\Sigma_M)$ is called \emph{trivial} if its range is reduced to $\{0\}$.
    \item A measure $\mu_M$ on a measurable space $(M,\Sigma_M)$ is called \emph{purely infinite} if its range is reduced to $\{0,\infty\}$.
    \item $\sigma(\mathcal{C})$ denotes the smallest $\sigma$-algebra that contains a family $\mathcal{C}$ of subsets of a set $M$.
    \item A subset of a measurable space $(M,\Sigma_M)$ is called \emph{measurable} when it belongs to $\Sigma_M$.
    \item $\mathcal{F}_{\mu_M}$ denotes the set of measurable subsets of a measure space $(M,\Sigma_M,\mu_M)$ with finite $\mu_M$-measure.
    \item $\Sigma_M|_B\coloneqq \left\{B\cap A: A\in \Sigma_M\right\}$ denotes the restriction of a $\sigma$-algebra $\Sigma_M$ on a set $M$ to subsets of $B\in \Sigma_M$.
    \item $\mu_M|_B$ denotes the restriction of a measure $\mu_M$ defined on a measurable space $(M,\Sigma_M)$ to $\Sigma_M|_B$ for some $B\in\Sigma_M$.
    \item $\mathcal{B}(M)$ denotes the \emph{Borel $\sigma$-algebra} on a topological space $M$.
    \item A measure $\mu_M$ on a topological space $M$ is called \emph{Borel} when it is defined on $\mathcal{B}(M)$.
    \item A subset $Z$ of a measure space $(M,\Sigma_M,\mu_M)$ is called \emph{$\mu_M$-null} if there exists $A\in \Sigma_M$ such that $Z\subset A$ and $\mu_M(A) = 0$. 
    \item $\mathcal{Z}_{\mu_M}$ denotes the set of all $\mu_M$-null sets of a measure space $(M,\Sigma_M,\mu_M)$.
    \item A property $\mathcal{P}$ that depends on the choice of point $x$ of a measure space $(M,\Sigma_M,\mu_M)$ is said to hold \emph{$\mu_M$-almost everywhere} ($\mu_M$-a.e. for short) if $\left\{x\in M: \text{$\mathcal{P}(x)$ is false}\right\}$ is a $\mu_M$-null set.
    \item In a measure space $(M,\Sigma_M,\mu_M)$, $\overline{\Sigma}_M$ denotes the \emph{completion of the $\sigma$-algebra $\Sigma_M$} with respect to $\mu_M$ and is defined as $\overline{\Sigma}_M\coloneqq \sigma(\Sigma_M \cup \mathcal{Z}_{\mu_M}) = \{ A\cup Z:(A,Z)\in \Sigma_M\times \mathcal{Z}_{\mu_M}\}$.
    \item In a measure space $(M,\Sigma_M,\mu_M)$, $\overline{\mu}_M$ denotes the \emph{completion of the measure $\mu_M$} and is defined as the mapping $\bar{\mu}_M: \overline{\Sigma}_M\to [0,\infty]$ such that $\bar{\mu}_M(A\cup Z)=\mu_M(A)$ for all $(A,Z)\in \Sigma_M\times \mathcal{Z}_{\mu_M}$. It is itself a measure on $(M,\overline{\Sigma}_M)$.
\end{itemize}
\label{sec:notations-end}

\section*{Acknowledgements}
\addcontentsline{toc}{section}{Acknowledgements}

This work is part of GS's PhD thesis. GS thus thanks his advisor Joan Alexis Glaun\`es (MAP5, Université Paris Cité) for his support, useful advice, and for introducing him to the subject of metamorphoses of manifold-valued images, which led him to dig into the theory of nonlinear Lebesgue spaces. The authors also thank Joan Alexis Glaun\`es for careful proofreading of early versions of this article as well as suggestions that helped improve the abstract and Th\'eo Dumont (LIGM, Université Gustave Eiffel) for suggestions that significantly helped improve the exposition of the final version of this article.

\appendix

\section{Omitted proofs}
\subsection{\texorpdfstring{Proof of \cref{prop:measurability_pointwise_limit}}{Proof of Proposition}}
\label{appendix:measurability_pointwise_limit}

\begin{proof}
    Let $(f_n)_{n\in\mathbb{N}}$ be a sequence in $\cLs(M,N)$ such that $f(x)\coloneqq \lim_{n\to \infty} f_n(x)$ exists in $N$ for all $x\in M$. The proof of the measurability of the pointwise limit is standard \cite[Proposition~8.1.10]{cohn2013measure}. For the separability of the range of $f$, observe that, since $f_n(M)$ is separable for all $n\in \mathbb{N}$, $\cup_{n\in\mathbb{N}}f_n(M)$ is separable as the countable union of separable sets and so is its closure $\overline{\cup_{n\in\mathbb{N}}f_n(M)}$. The set $f(M)$ is therefore separable as a subset of the separable set $\overline{\cup_{n\in\mathbb{N}}f_n(M)}$ \cite[(3.10.9)]{dieudonne1960treatise}.
\end{proof}

\subsection{\texorpdfstring{Proof of \cref{prop:equivalent_def_simple}}{Proof of Proposition}}
\label{appendix:equivalent_def_simple}

\begin{proof}
    Let $f: M\to N$ be a mapping.
    
    \par\medskip \noindent$(\Rightarrow)$ Assume $f\in \mathcal{E}(M,N)$. Then, there exists a finite indexing set $I$ such that $f(M)= \left\{ y_i \in N: i\in I\right\}$ and denote $M_i \coloneqq f^{-1}(\{y_i\})$. Thus, $(M_i)_{i\in I}$ forms a partition of $M$ and $M_i\in \Sigma_M$, by measurability of $f$.
    
    \par\medskip \noindent$(\Leftarrow)$ Assume there exists a finite partition $(M_i)_{i\in I}$ of $M$ such that $M_i\in \Sigma_M$ and $f|_{M_i} \equiv y_i \in N$. Then, $f(M) = \left\{ y_i: i\in I\right\}$ is finite and for every open set $U\subset N$ it holds that $f^{-1}(U)= \cup_{\{i\in I:y_i\in U\}}M_i\in \Sigma_M$, by stability of $\Sigma_M$ by finite union, that is, $f$ is measurable.
\end{proof}

\subsection{\texorpdfstring{Proof of \cref{prop:measurable_representative}}{Proof of Proposition}}
\label{appendix:proof_measurable_representative}
\begin{proof}
Let $f\in \overlinecL(M,N)$. 

\par\medskip \noindent$(\Rightarrow)$ Suppose that $f\in\overlinecLs(M,N)$. Since $f$ is $\mu_M$-essentially separably valued, there exists a measurable set $A\in \Sigma_M$ with $\mu_M(A)=0$ such that $f(M\setminus A)$ is separable. Fixing $z_0\in N$ and replacing $f$ by $z_0$ on $A$, we may assume that $f(M)$ is separable. Let $(y_k)_{k\geq 1}$ be a dense sequence in $f(M)$ and define for every $n\geq 1$, $f_n(x)=y_k$, where $k$ is the smallest index in $\{1,\ldots,n\}$ such that $d_N(f(x),y_k)=\min_{1\leq l\leq n}d_N(f(x),y_l)$. Equivalently, $f_n(x)=y_k$ if $d_N(f(x),y_k)<d_N(f(x),y_m)$ for $m<k$ and $d_N(f(x),y_k)\leq d_N(f(x),y_m)$ for $m>k$. It follows immediately that $f_n$ converges pointwise to $f$. Moreover,
for each $n\geq 1$, there exist a simple mapping $\tilde{f}_n\in\mathcal{E}(M,N)$ and a measurable set $A_n$ of null $\mu_M$-measure such that $\tilde{f}_n=f_n$ on $M\setminus A_n$. Consequently, $\tilde{f}_n$ converges pointwise to $f$ on $M\setminus A_\infty$, where $A_\infty\coloneqq\bigcup_n A_n$ has null $\mu_M$-measure. Define $\tilde{f}_\infty(x)=z_0$ on $A_\infty$ and $\tilde{f}_\infty(x)=\lim_{n\to \infty}\tilde{f}_n(x)$ elsewhere. Then, $\tilde{f}_\infty\in\mathcal{L}_s^0(M,N)$ (\cref{prop:measurability_pointwise_limit}) and $\tilde{f}_\infty\sim f$. 

\par\medskip \noindent$(\Leftarrow)$ Conversely, if there exists $\tilde{f}\in \overlinecLs(M,N)$ such that $\tilde{f}\sim f$, there is a $\mu_M$-null set $Z$ such that $\tilde{f}(M\setminus Z)$ is separable and $f(x)=\tilde{f}(x)$ for every $x\in M\setminus Z$. Hence, $f\in \overlinecLs(M,N)$.
\end{proof}

\subsection{\texorpdfstring{Proof of \cref{prop:countably_finite_radius}}{Proof of Proposition}}
\label{appendix:proof_countably_finite_radius}

\begin{proof}
  Let $f\in \cLs(M,N)$ and $\varepsilon >0$. Since $f(M)$ is separable, there exists a dense sequence $(y_n)_{n\in\mathbb{N}}$ in $f(M)$. Denote by $B(y,\varepsilon)$ the open ball of radius $\varepsilon$ for every $y\in N$. Define $g:M\to N$ by $g(x)\coloneqq y_{\tau(x)}$ for $x\in M$ where $\tau(x)\coloneqq\min\{\  n\geq 0,\ f(x)\in B(y_n,\varepsilon)\ \}$. For every $n\in \mathbb{N}$, $\tau^{-1}(\{n\})= f^{-1}(B(y_n,\epsilon))\cap(\bigcap_{k<n} f^{-1}(N\setminus B(y_k,\epsilon)))\in\Sigma_M$ where the intersection is understood to be $M$ when $n=0$. Hence, for every open set $U\subset N$, $g^{-1}(U) = \cup_{\{n\in\mathbb{N}: y_n \in U\}}\tau^{-1}(\{n\})\in \Sigma_M$. Therefore, $g$ is a countably valued measurable mapping. By construction, $D_\infty(g,f) \leq \varepsilon$.
\end{proof}
\subsection{\texorpdfstring{Proof of \cref{prop:indep_def}}{Proof of Proposition}}
\label{appendix:indep_def}

\begin{proof} \cref{prop:indep_def} is mentioned without proof in the literature, for instance in \cite[Section~3, p.~326]{sturm2002nonlinear}. 

\noindent Let $(h,h')\in \cLb(M,N)^2$ and $f\in \cLph(M,N)$. Then, by the triangle inequality, we have, for $\mu_M$-a.e. $x\in M$,
    $$d_N(f(x), h'(x)) \leq d_N(f(x), h(x)) + d_N(h(x),h'(x)).$$
For $p\in [1,\infty)$, we thus get, using the fact that $(a + b)^p \leq 2^{p-1}(a^p + b^p)$ for all $(a,b)\in \mathbb{R}_+^2$ and integrating over $M$, that
    $$D_p(f,h')^p \leq 2^{p-1}(D_p(f,h)^p + D_\infty(h,h')^p\mu_M(M)) < \infty.$$
Hence $L^p_h(M,N)\subset L^p_{h'}(M,N)$. For $p=\infty$, the same inclusion follows directly by taking the $\mu_M$-essential supremum over $M$. Interchanging $h$ and $h'$ yields the reverse inclusion. Therefore, $L^p_h(M,N) = L^p_{h'}(M,N)$.
\end{proof}

\subsection{\texorpdfstring{Proof of \cref{prop:trivial_lebesgue}}{Proof of Proposition}}
\label{appendix:proof_trivial_lebesgue}

\begin{proof}
Let $h\in \cLs(M,N)$.

\begin{enumerate}[label=(\roman*),wide]
    \item Let $p\in [1,\infty)$. 

\par\medskip \noindent$(\Rightarrow)$ We prove this implication by contraposition. Suppose that $\mu_M$ is purely infinite or $\lvert N\rvert=1$, and let $f\in\cLph(M,N)$. If $\lvert N\rvert=1$, then it is straightforward that $\Lph(M,N)=\{[h]\}$.
Assume now that $\mu_M$ is purely infinite. By Markov's inequality \cite[Proposition~2.3.10]{cohn2013measure}, for every $n\geq1$ and $A_n\coloneqq\{x\in M:d_N(f(x),h(x))\geq n^{-1}\}$, we have $n^{-p}\mu_M(A_n)\leq D_p(f,h)^p<\infty$. Hence, $\mu_M(A_n)<\infty$, and therefore $\mu_M(A_n)=0$ since $\mu_M$ is purely infinite. Since $A_\infty\coloneqq\{x\in M:d_N(f(x),h(x))>0\}=\bigcup_{n\geq1}A_n$, the continuity from below of the measure $\mu_M$ \cite[Proposition 1.2.5 (a)]{cohn2013measure} gives $\mu_M(A_\infty)=0$. Thus, $f\sim h$, and consequently $\Lph(M,N)=\{[h]\}$.

\par\medskip \noindent$(\Leftarrow)$ Suppose that $\mu_M$ is not purely infinite and that $\lvert N\rvert > 1$. Since $\mu_M$ is not purely infinite, there exists $A\in\Sigma_M$ such that $0<\mu_M(A)<\infty$. Since $\lvert N\rvert>1$, choose distinct points $(y_0,y_1)\in N^2$. For $n\geq0$, let $H_n\coloneqq\{x\in A\ :\ d_N(h(x),y_0)\leq n\}$. Since $A=\bigcup_{n\geq 1}H_n$, continuity from below of the measure $\mu_M$ \cite[Proposition 1.2.5 (a)]{cohn2013measure} gives $\mu_M(A)=\lim_{n\to\infty}\mu_M(H_n)$, and hence there exists $n_0>0$ such that $\mu_M(H_{n_0})>0$. For $i\in\{0,1\}$, set $A_i\coloneqq\{x\in H_{n_0}\ :\ h(x)\neq y_i\}$. Since $H_{n_0}=A_0\cup A_1$, there exists $i_0\in\{0,1\}$ such that $\mu_M(A_{i_0})>0$. Define $g:M\to N$ by $g=y_{i_0}$ on $H_{n_0}$ and $g=h$ on $M\setminus H_{n_0}$. Then $D_p(h,g)>0$, because $g\neq h$ on $A_{i_0}$, while $D_p(h,g)^p\leq\mu_M(H_{n_0})(n_0+d_N(y_0,y_{i_0}))^p<\infty$. Therefore, $[g]\neq[h]$ and $\Lph(M,N)\neq\{[h]\}$.
\item Let us prove each implication separately.

\par\medskip \noindent$(\Rightarrow)$ We prove this implication by contraposition. Suppose that $\mu_M$ is trivial or that $\lvert N\rvert=1$ and let $f\in \mathcal{L}^\infty_h(M,N)$.  When $\lvert N\rvert = 1$, it is straightforward that $L^\infty_h(M,N)=\{[h]\}$. When $\mu_M$ is trivial, $f\sim h$ for all $f\in \mathcal{L}^\infty_h(M,N)$ since $\mu_M(\{x\in M: f(x)\neq h(x)\}) \leq \mu_M(M) = 0$. Hence, $L^\infty_h(M,N)=\{[h]\}$.

\par\medskip \noindent$(\Leftarrow)$ Suppose that $\mu_M$ is nontrivial and that $\lvert N\vert> 1$. 
We can follow the construction used in case (i) in the simpler situation $A=M$. We still have $g\neq h$ on the set $A_{i_0}$, which has positive measure, and hence $D_\infty(g,h)>0$. Moreover, by the triangle inequality, $D_\infty(g,h)\leq n_0+d_N(y_0,y_{i_0})<\infty$. Therefore, $[g]\neq[h]$ and $L^\infty_h(M,N)\neq\{[h]\}$.\qedhere

\end{enumerate}

\end{proof}

\subsection{\texorpdfstring{Proof of \cref{prop:lebesgue_differ_from_base_bounded}}{Proof of Proposition}}
\label{appendix:proof_lebesgue_differ_from_base_bounded}

\begin{proof}
Let $p\in[1,\infty)$, $h\in\cLs(M,N)$, and $f\in\cLph(M,N)$. Fix $z_0\in N$ and, for every $n\in\mathbb{N}^*$, define the set $H_n\coloneqq \{x\in M\ :\ d_N(z_0,h(x))\leq n \text{ and } d_N(f(x),h(x))\geq n^{-1}\}$. Then $H_n\in\Sigma_M$, and Markov's inequality \cite[Proposition~2.3.10]{cohn2013measure} gives $n^{-p}\mu_M(H_n)\leq D_p(f,h)^p<\infty$, so that $\mu_M(H_n)<\infty$, and $h$ is bounded on $H_n$. Moreover, $\{x\in M\ :\ d_N(f(x),h(x))>0\ \}=\bigcup_{n\in\mathbb{N}^*}H_n$. Indeed, if $d_N(f(x),h(x))>0$, then, for $n$ sufficiently large, both $d_N(z_0,h(x))\leq n$ and $d_N(f(x),h(x))\geq n^{-1}$ hold. Hence, the set on which $f$ and $h$ differ is a countable union of measurable sets of finite $\mu_M$-measure on each of which $h$ is bounded.
    
\end{proof}

\subsection{\texorpdfstring{Proof of \cref{prop:continuity_restriction}}{Proof of Proposition}}
\label{appendix:continuity_restriction}

\begin{proof}
    Let $h\in \cLs(M,N)$, $p\in [1,\infty]$ and $(f,f')\in \cLph(M,N)^2$.
    \begin{enumerate}[label=(\roman*),wide]
        \item Let $B\in \Sigma_M$. If $p=\infty$, 
        \begin{align*}
        D_\infty(f|_B,f'|_B) &= \mu_M\text{-}\underset{x\in B}{\mathrm{ess}\,\sup}\, d_N(f|_B(x),f'|_B(x)) = \mu_M\text{-}\underset{x\in B}{\mathrm{ess}\,\sup}\, d_N(f(x),f'(x))\\
        &\leq \mu_M\text{-}\esssup_{x\in M} d_N(f(x),f'(x)) = D_\infty(f,f')
        \end{align*}
        and, if $p\in [1,\infty)$, 
        \begin{align*}
            D_p(f|_B,f'|_B) &= \left(\int_B d_N(f|_B(x),f'|_B(x))^p\dif\mu_M(x)\right)^{1/p}= \left(\int_B d_N(f(x),f'(x))^p\dif\mu_M(x)\right)^{1/p}\\
            &\leq \left(\int_M d_N(f(x),f'(x))^p\dif\mu_M(x)\right)^{1/p} = D_p(f,f').
        \end{align*}
        \item Similarly, we get, by the triangle inequality, that, if $p=\infty$,
        \begin{align*}
        \lVert \varphi_h(f) - \varphi_h(f')\rVert_{\infty,\mu_M} &= \mu_M\text{-}\esssup_{x\in M}\, \lvert d_N(f(x),h(x)) - d_N(f'(x),h(x))\rvert \\
        &\leq \mu_M\text{-}\esssup_{x\in M} d_N(f(x),f'(x)) = D_\infty(f,f')\
        \end{align*}
        and, if $p\in [1,\infty)$,
        \begin{align*}
            \lVert \varphi_h(f)-\varphi_h(f')\rVert_{p,\mu_M} &= \left(\int_M \lvert d_N(f(x),h(x)) - d_N(f'(x),h(x))\rvert ^p\dif\mu_M(x)\right)^{1/p}\\
            &\leq \left(\int_M  d_N(f(x),f'(x))^p\dif\mu_M(x)\right)^{1/p} = D_p(f,f').\qedhere
        \end{align*}
    \end{enumerate}
\end{proof}

\subsection{\texorpdfstring{Proof of \cref{lem:separation}}{Proof of Proposition}}
\label{appendix:separation}

\begin{proof}
\begin{enumerate}[label=(\roman*),wide]
    \item Let $n\geq 2$ be an integer. The case $n=2$ is normality. Assume it holds for $n-1$.  
Let $F=\bigcup_{i=1}^{n-1} C_i$. By normality, there exist disjoint opens $W_0$ and $U_n$ such that $F\subset W_0$ and $C_n\subset U_n$.  
By the induction assumption, there are disjoint opens $(V_i)_{1\leq i \leq n-1}$ satisfying $C_i\subset V_i$.  
Define $U_i:=V_i\cap W_0$ for all $1\leq i\le n-1$. Then, the $U_i$'s remain disjoint and contain the corresponding $C_i$'s, and lie inside $W_0$, hence are disjoint from $U_n$.  Thus $(U_i)_{1\leq i\leq n}$ is the desired family of open sets.

\item By normality, $C$ and $M\setminus U$ can be separated by disjoint opens $V$ and $W$. Then $\overline{V}\subset M\setminus W\subset U$, as desired.\qedhere
\end{enumerate}
\end{proof}

\subsection{\texorpdfstring{Proof of \cref{lem:separation_compact}}{Proof of Lemma}}
\label{appendix:separation_compact}

\begin{proof}
\begin{enumerate}[label=(\roman*),wide]
        \item The proof follows the same arguments as the proof of \ref{itm:separation_multiple} of \cref{lem:separation} by substituting normality with a result on the separation of compact sets in locally compact Hausdorff topological spaces \cite[Proposition 7.1.2]{cohn2013measure}.
        \item See \cite[Proposition 7.1.4]{cohn2013measure} for a proof.\qedhere
    \end{enumerate}
\end{proof}

\subsection{\texorpdfstring{Proof of \cref{cor:urysohn}}{Proof of Proposition}}
\label{appendix:ursyohn}

\begin{proof}
    Using \ref{itm:separation_compact} of \cref{lem:separation_compact}, we can find an open set $V$ with compact closure such that $K\subset V \subset \overline{V}\subset M\setminus F$. Using \ref{itm:separation_compact} of \cref{lem:separation_compact} one more time, we can find an open set $W$ with compact closure such that $K\subset W \subset \overline{W}\subset V$. Since $\overline{V}$ is compact and Hausdorff, it is normal (every closed subset being compact and by \ref{itm:separation_multiple_compact} of \cref{lem:separation_compact}). Therefore, we can apply \cref{th:urysohn} with the closed sets $ K$ and $\overline{V}\setminus W$, so that we get a continuous mapping $\tilde{I}: \overline{V} \to [0,1]$ such that $\tilde{I}|_K \equiv 1$ and $\tilde{I}|_{\overline{V}\setminus W}\equiv 0$. Then, define the mapping $I: M\to [0,1]$ such that $I|_{\overline{V}} = \tilde{I}$ and $I|_{M\setminus \overline{V}} \equiv 0$, which is continuous. Indeed, let $x\in M$. If $x\in V$, we already know that $I|_V = \tilde{I}|_V$. If $x\in M\setminus V$, then $x\in M\setminus \overline{W}$, so that there exists a neighborhood $U_x$ of $x$ such that $U_x \subset M\setminus \overline{W}$ and $I|_{U_x} \equiv 0$. In both cases, there exists a neighborhood $U_x$ of $x$ such that $I|_{U_x}$ is continuous. Since $F\subset M\setminus\overline{V}$, we also have $I|_F\equiv 0$.
\end{proof}

\section{Omitted results}

\renewcommand{\theproposition}{\Alph{section}.\arabic{proposition}} 

\subsection{\texorpdfstring{Explicit expression of the dense subspace in \cref{th:density_simple}}{Explicit expression of the dense subspace in Theorem}}

When the base mapping is simple, the set of simple mappings that are Lebesgue is the set of almost simple mappings.

\begin{proposition}
\label{prop:simple_lebesgue_are_alsmot_simple}
    Let $h\in E(M,N)$ and $p\in [1,\infty)$. Then, $E(M,N)\cap \Lph(M,N) = E_h(M,N)$.
\end{proposition}
\begin{proof}
Let $h\in \mathcal{E}(M,N)$ and $p\in [1,\infty)$. The inclusion $E_h(M,N)\subset E(M,N)\cap \Lph(M,N)$ is obvious since $h$ is simple and almost simple mappings differ from $h$ on a set of finite $\mu_M$-measure. Let us prove the reverse inclusion: Let $g\in\mathcal{E}(M,N)\cap\cLph(M,N)$. Then there exist finite indexing sets $I$ and $J$ such that $g(M)=\{y_i\in N:i\in I\}$ and $h(M)=\{z_j \in N: j\in J\}$. For every $i\in I$ and $j\in J$, set $G_i\coloneqq g^{-1}(\{y_i\})$ and $H_j\coloneqq h^{-1}(\{z_j\})$. For every $(i,j)\in I\times J$ such that $y_i\neq z_j$, we have $d_N(y_i,z_j)^p\mu_M(G_i\cap H_j)\leq D_p(g,h)^p<\infty$, and hence $\mu_M(G_i\cap H_j)<\infty$. Since $\{x\in M:g(x)\neq h(x)\}$ is the union of the $G_i\cap H_j$'s such that $y_i\neq z_j$, this set is a finite union of measurable sets of finite measure. Therefore, $g\in\mathcal{E}_h(M,N)$.
\end{proof}

\section{Reminders}

\subsection{Measure theory}

The Borel $\sigma$-algebra of a separable metric space is countably generated.

\label{appendix:borel_separable}

\begin{proposition}
\label{prop:borel_separable}
    Let $(M,d_M)$ be a separable metric space. Then, $\mathcal{B}(M)$ is countably generated.
\end{proposition}
\begin{proof}
    Since $(M,d_M)$ is a separable metric space, it is second-countable \cite[Appendix D.32]{cohn2013measure}. Thus, let $(U_n)_{n\in\mathbb{N}}$ be a countable basis of $M$, so that for every open subset $V\subset M$ there exists a set of indices $I\subset \mathbb{N}$ such that $V = \cup_{i\in I} U_i$. Therefore, if $\mathcal{F} \coloneqq \sigma(\left\{U_n:n\in \mathbb{N}\right\})$ denotes the $\sigma$-algebra generated by the countable basis of $M$, the basis property yields that $V\in \mathcal{F}$ and thus that $\left\{V\subset M: \text{$V$ is open}\right\} \subset \mathcal{F}$. Hence, $\mathcal{B}(M)$ being the smallest $\sigma$-algebra containing the open subsets of $M$, we get that $\mathcal{B}(M)\subset \mathcal{F}$. The reverse inclusion follows from observing that $\left\{U_n:n\in \mathbb{N}\right\} \subset \mathcal{B}(M)$ and the fact that $\mathcal{F}$ is the smallest $\sigma$-algebra containing $\left\{U_n:n\in \mathbb{N}\right\}$. Thus, $\mathcal{B}(M)$ is countably generated.
\end{proof}

\subsection{Differentiable manifolds}

\renewcommand{\thedefinition}{\Alph{section}.\arabic{definition}} 
\renewcommand{\thetheorem}{\Alph{section}.\arabic{theorem}} 
\renewcommand{\theremark}{\Alph{section}.\arabic{remark}} 

Let us first recall a definition of \emph{atlases} by paraphrasing \cite[Chapter~II,\S 1]{lang2012fundamentals}.

\begin{definition}[Atlases]
\label{def:atlases}
Let $r\in \mathbb{N}^*\cup \{\infty\}$. An \emph{atlas of class $\mathcal{C}^r$} on a set $M$ is a family of pairs $\{(U_\alpha,\varphi_\alpha)\}_{\alpha\in A}$ indexed by a set $A$ satisfying the following conditions:
\begin{enumerate}[label=(\roman*)]
    \item each $U_\alpha$ is a subset of $M$ and the $U_\alpha$'s cover $M$.
    \item each $\varphi_\alpha$ is a bijection from $U_\alpha$ onto an open subset $\varphi_\alpha (U_\alpha)$ of some Banach space $\mathbb{B}_\alpha$ and for all $(\alpha,\alpha')\in A^2$ the set $\varphi_\alpha(U_\alpha \cap U_{\alpha'}) $ is open in $\mathbb{B}_\alpha$.
    \item the mapping $\varphi_\alpha \circ \varphi_{\alpha'}^{-1}: \varphi_{\alpha'}(U_\alpha\cap U_{\alpha'}) \to \varphi_\alpha(U_\alpha\cap U_{\alpha'})$ is a $\mathcal{C}^r$ isomorphism for each pair of indices $(\alpha,\alpha')\in A^2$.
\end{enumerate}
Two atlases $\{(V_\beta, \psi_\beta)\}_{\beta\in B}$ and $\{(U_\alpha,\varphi_\alpha)\}_{\alpha\in A}$ on $M$ indexed by some sets $A$ and $B$ are said to be \emph{compatible} if each pair $(V_\beta, \psi_\beta)$ is such that the mapping $\varphi_\alpha \circ \psi_\beta^{-1}: \psi_\beta(U_\alpha\cap V_\beta)\to \varphi_\alpha(U_\alpha\cap V_\beta)$ is a $\mathcal{C}^r$ isomorphism for all pairs $(U_\alpha,\varphi_\alpha)$ and conversely. This relation defines an equivalence relation between atlases on $M$.
\end{definition}
\begin{remark}[On topologies induced by atlases]
    Note that an atlas uniquely defines a topology in which the charts $\{(U_\alpha,\varphi_\alpha)\}_{\alpha\in A}$ are such that $U_\alpha$'s are open sets and the $\varphi_\alpha$'s are topological isomorphisms. This induced topology is not required to be Hausdorff.
\end{remark}
 Using the equivalence relation of compatibility between atlases, \emph{Banach manifolds} are then defined as follows:

\begin{definition}[Banach manifolds]
\label{def:banach_manifolds}
    Let $r\in \mathbb{N}^*\cup \{\infty\}$. A \emph{$\mathcal{C}^r$ Banach manifold} is a set $M$ together with an equivalence class of atlases of class $\mathcal{C}^r$ (\cref{def:atlases}) on this set.
\end{definition}

One can then define smooth mappings between two Banach manifolds using charts.

\begin{definition}[Smooth mappings]
\label{def:smooth_mappings}
    Let $r\in \mathbb{N}^*\cup \{\infty\}$. Given two $\mathcal{C}^r$ Banach manifolds $M$ and $N$, a mapping $f: M\to N$ is said to be of class $\mathcal{C}^r$ if for every point $x\in M$ there exists a chart $(U,\varphi)$ such that $x\in U$ and a chart $(V,\psi)$ such that $f(x)\in V$ and $f(U)\subset V$ for which $\psi \circ f \circ \varphi^{-1}: \varphi(U)\to \psi(V)$ is of class $\mathcal{C}^r$ in the usual sense \cite[Chapter~I,\S 3]{lang2012fundamentals}. The set of such mappings is denoted by $\mathcal{C}^r(M,N)$. 
\end{definition}

On Banach manifolds, one is often interested in using \emph{partitions of unity} to define mappings locally and extend them to the whole manifold. We thus recall a definition inspired by the exposition in \cite[Chapter~5, p.~300]{engelking1977general}.

\begin{definition}[Partitions of unity]
\label{def:partition_unity}
    Let $r\in \mathbb{N}^*\cup \{\infty\}$. Suppose $M$ is a $\mathcal{C}^r$ Banach manifold (not necessarily Hausdorff). Then,
    \begin{enumerate}[label=(\roman*)]
        \item a family $\{\rho_s\}_{s\in S}$ indexed by some set $S$ of $\mathcal{C}^r$ mappings from $M$ to $[0,1]$ is called a \emph{partition of unity} if $\sum_{s\in S} \rho_s(x) = 1$ for all $x\in M$ and we denote by $\operatorname{carr}(\rho_s)\coloneqq \{x\in M: \rho_s(x) \neq 0\}$ the \emph{carrier} of $\rho_s$.
        \item a partition of unity $\{\rho_s\}_{s\in S}$ is called \emph{locally finite} if for all $x\in M$ there exist a neighborhood $U_x\subset M$ and  a finite set $S_x \subset S$ such that for all $x'\in U_x$ we have $\rho_s(x')=0$ for all $s\in S\setminus S_x$ and $\sum_{s\in S_x} \rho_s(x')=1$. This is equivalent to saying that the cover $\{\operatorname{carr}(\rho_s)\}_{s\in S}$ of $M$ is \emph{locally finite}, that is, for all $x\in M$ there exists a neighborhood $U_x\subset M$ such that $\{s\in S:\operatorname{carr}(\rho_s)\cap U_x \neq \emptyset\}$ is finite.
        \item a partition of unity $\{\rho_s\}_{s\in S}$ is called \emph{subordinated} to an open cover $\{U_\alpha\}_{\alpha\in A}$ of $M$ indexed by some set $A$ if $\{\operatorname{carr}(\rho_s)\}_{s\in S}$ is a \emph{refinement} of $\{U_\alpha\}_{\alpha\in A}$, that is, for all $s \in S$ there exists $\alpha \in A$ such that $\operatorname{carr}(\rho_s)\subset U_\alpha$.
    \end{enumerate}
\end{definition}

\bibliographystyle{alpha}
\bibliography{bibfile_final}

\end{document}